\documentclass[twoside,leqno,fleqn]{article}
\usepackage{amsmath}
\usepackage{amssymb}
\usepackage{array}
\usepackage[shortlabels]{enumitem}
\usepackage[onehalfspacing]{setspace}
\usepackage[T1]{fontenc}
\usepackage{graphicx}
\usepackage{xcolor}
\usepackage{hyperref}
 
\newcommand{\bbbn}{\mathbb{N}} 
\newcommand{\qed}{\hfill{$\Box$}} 
\newtheorem {overview}{Overview} 
\newtheorem{satz}{Satz}[section]

\newtheorem{corollary}[satz]{Corollary}
\newtheorem{theorem}[satz]{Theorem}
\newtheorem{remark}[satz]{Remark}
\newtheorem{example}[satz]{Example}

\counterwithout{equation}{section}
\begin{document}
\pagestyle{headings} 
\thispagestyle{empty}

\title{Abstract computation over first-order structures. Part II\,a: Moschovakis' operator and other non-determinisms}
\author{Christine Ga\ss ner}

\begin{center}

{\Large Abstract Computation over First-Order Structures \vspace{0.2cm}}\\{\Large Part II\,a\\ Moschovakis' Operator and Other Non-Determinisms}
\vspace{0.7cm}\\
{\bf Christine Ga\ss ner} \footnote{In Part III and Part IV, we present a generalization of the results discussed at the CCC in Kochel in 2015. My thanks go to Dieter Spreen, Ulrich Berger, and the other organizers of the conference in Kochel. Further results were discussed at the CCA 2015. I would like to thank the organizers of this conference and in particular Martin Ziegler and Akitoshi Kawamura.} {\vspace{0.1cm}\\University of Greifswald, Germany, 2025\\ gassnerc@uni-greifswald.de}\\\end{center} 

\vspace{0.3cm}

\begin{abstract} 

 BSS RAMs over first-order structures help to characterize algorithms for processing objects by means of useful operations and relations. They are the result of a generalization of several types of abstract machines. We want to discuss whether this concept that allows a machine-oriented characterization of algorithms is sufficiently general for describing also other models of computation. Yiannis N. Moschovakis introduced a concept of abstract computability of functions on the basis of recursive definability over first-order structures. Moschovakis' search operator is the counterpart to the operator introduced by Stephen C. Kleene and suitable for structures without computable minima. To compare our concept with Moschovakis' generalization of the theory of recursive functions, we extend the abilities of BSS RAMs by an operator that makes it possible to provide information about computable functions and their inverses in a non-deterministic way. In Part II\,b, we compare several non-determinisms, summarize effects resulting from the restriction of guesses to constants, and take into account properties such as the semi-decidability of oracle sets, the semi-decidability of the identity relation, and the recognizability of constants. 

\vspace{0.4cm}

\end{abstract}

 \vspace*{0.3cm}

\newpage

\tableofcontents

\vspace{0.5cm}

\noindent {\large\bf Contents of Part I} 

{\small
\begin{description}[itemindent=-1cm]\parskip -0.6mm

\item {\bf \hspace{0.1cm}1\quad Introduction} 

\item {\bf \hspace{0.1cm}2\quad BSS-RAMs over first-order structures} 

First-order structures and their signatures. The $\sigma$-programs. Finite machines and more. Infinite machines and their configurations. Computation and semi-decidability by BSS RAMs. Non-deterministic BSS RAMs. Oracle machines. Other non-deterministic BSS RAMs.

\item {\bf \hspace{0.1cm}3\quad Examples} 

Input and output functions. Computing functions and deciding sets. Evaluating formulas. 
\end{description}}

\noindent {\large\bf Contents of Part II\,b} 

{\small
\begin{description}[itemindent=-1cm]\parskip -0.6mm
\item {\bf \hspace{0.1cm}6\quad Two constants and consequences}

Two constants as outputs and the decidability. Two constants and their recognizability. Binary and other non-determinisms restricted to constants. Different non-determinisms and relationships. 
\end{description}
\noindent For Part I (including a summary and an outlook), see \cite{GASS25A}. Part II\,b (including an introduction and repetition, a summary, and an outlook) is coming soon.}

\newpage
\section*{Introduction and repetition} 
\addcontentsline{toc}{section}{\bf Introduction and repetition}
\markboth{INTRODUCTION AND REPETITION}{INTRODUCTION AND REPETITION}
On the one hand, we try to characterize the classical term {\sf algorithm} in a broad and comprehensive sense. By using familiar terms from computer science, we describe the abstract computability over any first-order structure in a machine-oriented way even when this does not mean that the algorithms can be implemented for executing them by machines (today). On the other hand, we try to use only as much types of instructions as necessary and to use an easily manageable memory space with a small but sufficient capacity. The advantage can be that we only have to consider a few cases if we want to prove the undecidability of a decision problem, and so on.

\vspace{0.1cm}

{\bf Types of instructions.} The basic types of instructions are summarized in Overview \ref {Sigma_InstructionsInfini}. $\sigma$ can be a signature of the form $(n_1;m_1, \dots,m_{n_2}; k_1, \dots,k_{n_3})$. 
\begin{overview}[$\sigma$-instructions and other instructions]\label{Sigma_InstructionsInfini}

\hfill

\nopagebreak 

\noindent \fbox{\parbox{11.8cm}{

\vspace{0.1cm}

Computation instructions

\hspace{0.6cm} $(1)$ \quad $\ell : \, Z_j:= f_i^{m_i}(Z_{j_1},\ldots, Z_{j_{m_i}}) $

\hspace{0.6cm} $(2)$ \quad $\ell : \,Z_j:=c_i^0$

\vspace{0.1cm}

Copy instructions

\hspace{0.6cm} $(3)$ \quad $\ell : \,Z_{I_j}:=Z_{I_k}$

\vspace{0.1cm}

Branching instructions

\hspace{0.6cm} $(4)$ \quad {\sf $\ell : \,$ if $r_i^{k_i}(Z_{j_1},\ldots, Z_{j_{k_i}})$ then goto $\ell _1$ else goto $\ell _2$}

\vspace{0.1cm}

Index instructions

\hspace{0.6cm} $(5)$ \quad {\sf $\ell : \,$ if $I_{j}=I_k$ then goto $\ell _1$ else goto $\ell _2$} 

\hspace{0.6cm} $(6)$ \quad $\ell : \,I_j:=1$

\hspace{0.6cm} $(7)$ \quad $\ell : \,I_j:=I_j+1$

\vspace{0.1cm}

Stop instruction

\hspace{0.63cm} $(8)$ \quad$l : $ {\sf stop}
}}
\end{overview}
The interpretation of each $\sigma$-instruction can be done by means of a first-order structure ${\cal A}$ if ${\cal A}$ includes a reduct of signature $\sigma$, or in other words, if ${\cal A}$ is an expansion of a structure of signature $\sigma$. For ${\cal A}$-machines, see Part I \cite[p.\,9]{GASS25A}.

\begin{overview}
[Oracle instructions]

\hfill

\nopagebreak 

\noindent \fbox{\parbox{ 11.8cm}{

\vspace{0.1cm}

Oracle instructions with queries

\hspace{0.6cm} $(9)$ \quad {\sf $\ell :\,$ if $(Z_1,\ldots,Z_{I_1})\in \!{ \cal O}$ then goto $\ell_1$ else goto $\ell_2$}

\vspace{0.1cm}

Oracle instructions with $\nu$-operator

 \hspace{0.43cm} $(10)$ \quad $\ell:\,\, Z_j:=\nu[{\cal O}](Z_1,\ldots,Z_{I_1})$
}}
\end{overview}
${\cal O}$ stands for any {\sf oracle set} $Q\subseteq U_{\cal A}^\infty$. The evaluation of queries of the form {\em $(c(Z_1),\ldots, c(Z_{c(I_1)})) \in Q$?} is possible by using instructions of type (9). The execution of an algorithm can therefore depend on states belonging to an external space of states which we call {\em oracle set} and which can also be considered as a {\em controlled space of search}. We introduced them in \cite{GASS97} and considered them in \cite{GASS07, GASS08A, GASS08B, GASS08C, GASS08D, GASS09A, GASS09B, GASS10} and in Part I \cite{GASS25A}. Instructions of type (10) enable to select components of a tuple belonging to this controlled space of search. We investigate them in more detail in the following sections.

\vspace{0.2cm}

{\bf Classes of structures.}
Let ${\sf Struc}$ be the class of all first-order structures of any signature. Let {\em $c_1$ and $c_2$} be any two different constants, e.g., $1$ {\sf and} $2$ or $a$ {\sf and} $b$, etc. For such constants, let ${\sf Struc}_{c_1,c_2}$ be the subclass of ${\sf Struc}$ containing all structures $(U; (c_i)_{i\in N_1}; (f_i)_{i\in N_2}; (r_i)_{i\in N_3}) $ with $\{1,2\} \subseteq N_1$ (which means that $\{c_1,c_2\} \subseteq \{c_i\mid i\in N_1\}$ holds). Thus, the formula ${\cal A}\in {\sf Struc}_{c_1,c_2}$ expresses that ${\cal A}$ is a structure containing at least two constants denoted by $c_1$ and $c_2$.

\vspace{0.2cm}

{\bf $\sigma$-Programs.} Any so-called program ${\cal P}$ of an ${\cal A}$-machine has the form

\vspace{0.1cm}

$\!\!1: {\sf instruction}_{1}; \quad \ldots\,;\quad\ell _{\cal P}-1:\, {\sf instruction}_{\ell _{\cal P}-1};\quad \,\,\,\ell _{\cal P}:\, {\sf stop}.\hfill \textcolor{blue}{(*)}\,\,$\label{star}

\vspace{0.1cm}

\noindent It is a finite sequence of  \big[$\sigma$-\big]instructions and their labels. The execution of the program ${\cal P}_{\cal M}$ by an ${\cal A}$-machine ${\cal M}$ is, by definition, the stepwise transformation of so-called configurations by means of the transition system $({\sf S}_{\cal M},\to_{\cal M} )$. 

\vspace{0.2cm}

{\bf Configurations.} Each configuration of an ${\cal A}$-machine ${\cal M}$ describes the overall state of the machine at a certain point of time. Any configuration  of an infinite 1-tape ${\cal A}$-machine ${\cal M}$ is a sequence $(\ell, \nu_1,\ldots,\nu_{k_{\cal M}}, u_1,u_2,\ldots)$ consisting of an instruction counter and states stored in $k_{\cal M}$ index registers and in an infinite number of $Z$-registers and can be denoted by $(\ell\,.\,\vec \nu\,.\,\bar u)$ where $\ell \in \{1,\ldots, \ell_{{\cal P}_{\cal M}}\}$, $\vec \nu=( \nu_1,\ldots,\nu_{k_{\cal M}})\in \mathbb{N}_+^{k_{\cal M}}$, and $\bar u=( u_1,u_2,\ldots)\in U_{\cal A}^\omega$ hold. The transformation of one configuration $con_t$ into the next configuration $con_{t+1}$ is the change of a finite number of components resulting from the application of certain functions introduced for describing the execution of instructions purely mathematically. For the transformation of configurations, see also \cite[Overviews 20 and 21]{GASS25A}. 

\vspace{0.2cm}

{\bf Simulation of an ${\cal A}$-machine.} We speak of the {\sf simulation of an ${\cal A}$-machine ${\cal M}$ on $\vec x$ by a machine ${\cal N}$} (over any structure) when the computation by ${\cal N}$ on a suitable input $f(\vec x)$ is a computation process that allows to compute (the codes of) all individuals and indices that are important intermediate results for getting the expected states or results of ${\cal M}$. These can be a label such as $\ell_{{\cal P}_{\cal M}}$ or a result such as ${\rm Res}_{\cal M}(\vec x)$ computed by ${\cal M}$ so that ${\rm Res}_{\cal M}(\vec x)=g^{-1}({\rm Res}_{\cal N}(f(\vec x)))$ holds for suitable $f$ and $g$ and a function ${\rm Res}_{\cal N}$ defined e.g.\,\,as in \cite[p.\,585]{GASS20}. In this sense, each multi-tape BSS RAM can be simulated by a 1-tape machine so that both machines have the same result function (without the use of additional homomorphisms or other mappings such as $f$ or $g$). More precisely, each $d$-tape machine can be simulated by a 1-tape machine working in a multi-tape mode with the help of $d$ tracks (cf.\,\,\cite{GASS20}). Therefore, we generally consider only 1-tape machines. However, to describe the simulation of a 1-tape machine by a universal machine or a machine of another type, we for simplicity use a multi-tape machine when a description of this kind is easier to follow, since we know that such a machine can also be simulated by a 1-tape machine again.

\vspace{0.2cm}

{\bf Multi-tape machines.} In \cite{GASS20}, we introduced deterministic $d$-tape machines for storing individuals on several tapes given by $Z_{d_i,1}, Z_{d_i,2},\ldots $ for $d_i\leq d$. 

For any $d\geq 1$, the non-deterministic $d$-tape machines can be analogously introduced as the non-deterministic 3-tape BSS RAMs in the following overview.

\begin{overview}
[Non-deterministic 3-tape BSS RAMs in $({\sf M}^{(3)}_{\cal A})^{\rm ND}$]\label{NDThreeTape}

\hfill

\nopagebreak 

\noindent \fbox{\parbox{11.8cm}{
A {\em $3$-tape ${\cal A}$-machine} ${\cal M}$ is a tuple $((U_{\cal A}^\omega)^3,\mathbb{N}_+^{\kappa_1}\times\mathbb{N}_+^{\kappa_2}\times\mathbb{N}_+^{\kappa_3}, {\cal L}, {\cal P}, {\cal B}, {\rm In}, {\rm Out})$ equipped, for all $d\leq 3$, with 
\begin{itemize}
\item\parskip -0.5mm $Z$-registers $Z_{d,1},Z_{d,2},\ldots$ forming the $d^{\rm th}$ tape,
\item $\kappa_d$ index registers $I_{d,1},\ldots,I_{d,\kappa_d}$.
\end{itemize}

A pair $(\ell,(\vec \nu^{\,(d)}\,.\,\bar u^{(d)})_{d=1..3})$ is a {\em configuration} of ${\cal M}$ if 

\begin{itemize} 
\item\parskip -0.5mm 
$\ell$ is the label $c(B)$ stored in the register $B$ of ${\cal M}$, 
\item the component $\nu_{d,k}$ of $\vec \nu^{\,(d)}$ is the content $c(I_{d,k})$ of $I_{d,k}$ in ${\cal M}$ ($k\leq \kappa_d$), 
\item the component $u_{d,k}$ of $\bar u^{(d)}$ is the content $c(Z_{d,k})$ of $Z_{d,k}$ in ${\cal M}$. 
\end{itemize}

A $3$-tape ${\cal A}$-machine ${\cal M}$ is a {\em non-deterministic $3$-tape BSS RAM} if 
\begin{itemize} 
\item\parskip -0.5mm 
each $(\vec x,(\vec \nu^{\,(d)}\,.\,\bar u^{(d)})_{d=1..3})\in {\rm In}_{\cal M}$ is a pair determined by
\begin{itemize} 
\item\parskip -0.5mm 
$ \vec x= (x_1, \ldots, x_n)\in U_{\cal A}^\infty$,
\item $\vec \nu^{\,(1)} = (n,1,\ldots, 1)\in\mathbb{N}_+^{\kappa_1}$, 
\item $\bar u^{(1)} \,= (x_1, \ldots ,x_n, y_1,\ldots, y_m,x_n,x_n,\ldots )$ \hfill with $(y_1,\ldots, y_m) \in U_{\cal A}^\infty$, 
\item $\vec \nu^{\,(2)}=(1,\ldots , 1) \in\mathbb{N}_+^{\kappa_2}$ and $\bar u^{(2)} =(x_n,x_n,\ldots )$, 
\item $\vec \nu^{\,(3)}=(1,\ldots , 1) \in\mathbb{N}_+^{\kappa_3}$ and $\bar u^{(3)} =(x_n,x_n,\ldots )$,
\end{itemize}
\item ${\rm Out}_{\cal M}$ 
is defined by ${\rm Out}_{\cal M}((\vec \nu^{\,(d)}\,.\,\bar u^{(d)})_{d=1..3})= (u_{1,1}, \ldots,u_{1,\nu_{1,1}})$. 
\end{itemize}Let $({\sf M}^{(3)}_{\cal A})^{\rm ND}$ be the class of all non-deterministic $3$-tape BSS RAMs. 
}}
\end{overview}

\vspace{0.1cm}

\begin{overview}[Instructions for $d$-tape machines]\label{Sigma_Instructionsd_tape}

\hfill

\nopagebreak 

\noindent \fbox{\parbox{11.8cm}{

\vspace{0.1cm}

Computation instructions \hfill $^{(d_0,d_1,d_2,\ldots \leq d)}$

\hspace{0.6cm}$^{(d)}(1)$ \quad $\ell : \, Z_{d_0,j}:= f_i^{m_i}(Z_{d_1,j_1},\ldots, Z_{d_{m_i},j_{m_i}})$

\hspace{0.6cm}$^{(d)}(2)$ \quad $\ell : \,Z_{d_0,j}:=c_i^0$

\vspace{0.2cm}

 Copy instructions

\hspace{0.6cm}$^{(d)}(3)$ \quad $\ell : \,Z_{d_1,I_{d_1,j}}:=Z_{d_2,I_{d_2,k}}$

\vspace{0.2cm}

Branching instructions

\hspace{0.6cm}$^{(d)}(4)$ \quad {\sf $\ell : \,$ if $r_i^{k_i}(Z_{d_1,j_1},\ldots, Z_{d_{k_i},j_{k_i}})$ then goto $\ell _1$ else goto $\ell _2$}\hspace*{0.5cm}

\vspace{0.2cm}

Index instructions

\hspace{0.6cm}$^{(d)}(5)$ \quad {\sf $\ell : \,$ if $I_{d_1,j}=I_{d_2,k}$ then goto $\ell _1$ else goto $\ell _2$} 

\hspace{0.6cm}$^{(d)}(6)$ \quad $\ell : \,I_{d_1,j}:=1$

\hspace{0.6cm}$^{(d)}(7)$ \quad $\ell : \,I_{d_1,j}:=I_{d_1,j}+1$

\vspace{0.2cm}

Stop instruction

\hspace{0.6cm}$^{(d)}(8)$ \quad $l : $ {\sf stop}

\vspace{0.2cm}

Oracle instructions with queries

\hspace{0.6cm}$^{(d)}(9)$ \quad {\sf $\ell :\,$ if $(Z_{d_1,1},\ldots,Z_{d_1,I_{d_1,1}})\in { \cal O}$ then goto $\ell_1$ else goto $\ell_2$}

\vspace{0.2cm}

Oracle instructions with $\nu$-operator

\hspace{0.42cm}$^{(d)}(10)$ \quad $\ell:\, Z_{d_1,j}:=\nu[{\cal O}](Z_{d_2,1},\ldots,Z_{d_2,I_{d_2,1}})$

\vspace{0.1cm}

}}
\end{overview}

\vspace{0.1cm}

In Section \ref{SectBSS_Moschov}, we give more background information that can help to better understand what the benefit of the $\nu$-instructions with Moschovakis' operator is. The section contains examples, that demonstrate the proper extension of a class of computable functions by the use of the $\nu$-operator for a first-order structure, and parts of the proofs of the main theorems in a less formal style. For semi-decidable oracle sets to which the $\nu$-operator is applied and weaker variants of oracles, we discuss the possibilities to simulate the execution of this kind of instructions by means of non-deterministic BSS RAMs. More details of the proofs are presented in Section \ref{SectSimMachinesWithNuInst}. Finally, we give a summary and a small outlook. 

In Section 6 (Part II\,b), we discuss the power of non-deterministic machines that may only use constants as guesses and compare the binary non-determinism, the digital non-determinism, and other non-determinisms for several classes of structures. 

We are doing all these to better understand {\sf the term algorithm}. 

\setcounter{section}{3}

\section{BSS RAMs with Moschovakis' operator}\label{Section3MoschoOper}\label{SectBSS_Moschov}

The guessing process executed by a non-deterministic BSS RAM can be simulated by means of a non-de\-terministic operator. Yiannis N. Moschovakis introduced a search operator $\nu$ in \cite[p.\,\,449]{MOSCHO}.\footnote{Therefore, we will also use this notation. There is no connection between the operator $\nu$ and the content $\nu_1$ of $I_1$, and so on.} This $\nu$-operator can be regarded as a variant derived from Kleene's $\nu$-operator (see \cite[p.\,347]{KLEENE}) and it can be used instead of Kleene's $\mu$-operator introduced for providing a minimum of certain subsets of $\bbbn$ (see \cite{KLEENE}). If a structure does not allow to determine the minimum for any subset of the universe, we should not use an operator such as $\mu$ with properties described by Stephen C. Kleene. Therefore, we follow Moschovakis' idea and use a non-deterministic operator here also denoted by $\nu$. For any structure ${\cal A}\in {\sf Struc}$ and any partial function $f: \,\, \subseteq\! U_{\cal A}^\infty\to U_{\cal A}^\infty$, let us first consider

\vspace{0.2cm}

$\nu[f](\vec z)=_{\rm df}\{ y_1\in U_{\cal A}\mid f(\vec z\,.\,y_1) \downarrow \mbox{ or }(\exists (y_2,\ldots,y_m)\in U_{\cal A}^\infty)( f(\vec z\,.\, \vec y) \downarrow)\}$

\vspace{0.2cm}

\noindent for any $\vec z\in U_{\cal A}^\infty$ where $f(\vec u)\downarrow$ means that $f(\vec u)$ is defined and $\exists (y_2,\ldots,y_m)\in U_{\cal A}^\infty$ stands for $(\exists m\geq 2) (\exists (y_2,\ldots,y_m)\in U_{\cal A}^m)$ where we use $\vec y=(y_1,\ldots,y_m)$ and $(\vec z\,.\, \vec y)=(z_1,\ldots,z_{n_0},y_1,\ldots,y_m)$ for $\vec z=(z_1,\ldots,z_{n_0})$. A suitable non-deterministic oracle instruction with Moschovakis' operator could be 

\vspace{0.2cm}

$\ell:\,\, Z_j:=\nu[{\sf f}](Z_1,\ldots,Z_{I_1})$

\vspace{0.2cm}

\noindent where ${\sf f}$ stands for a partial function $f:\,\,\subseteq\! U_{\cal A}^\infty \to U_{\cal A}^\infty $ computable over ${\cal A}$. Then, its execution by an ${\cal A}$-machine ${\cal M}$ has the consequence that $ Z_j$ gets a value that belongs to $ {\nu}[f](c(Z_1),\ldots,c(Z_{n_0}))$ if $c(I_1)=n_0$ and $ {\nu}[f](c(Z_1),\ldots,c(Z_{n_0}))\not= \emptyset$. Otherwise, this instruction has the effect that ${\cal M}$ loops forever (which means that ${\cal M}$ does not halt). A similar type of instructions was introduced for ${\cal A}$-machines that can use constants $c_1, c_2,\ldots$ of ${\cal A}$. In \cite{GaVV15}, $\nu[f]$ was defined by 

\vspace{0.2cm}

$\nu[f](\vec z)=\{ y_1\in U_{\cal A}\mid f(\vec z\,.\,y_1)=c_1 \mbox{ or }(\exists (y_2,\ldots,y_m)\in U_{\cal A}^\infty)( f(\vec z\,.\, \vec y) =c_1)\}$.

\vspace{0.2cm}

\noindent 
Because $f$ can also be the partial characteristic function $\bar \chi_Q: \, \subseteq U_{\cal A}^\infty\to \{c_1\}$ of a set $Q\subseteq U_{\cal A}^\infty$ (with $\bar \chi_Q(\vec z)=c_1$ if $\vec z \in Q$ and $\bar \chi_Q(\vec z)\uparrow$ otherwise), it is also possible to extend the applicability of the $\nu$-operator. We use the definition as given in \cite{GASS25A} and discuss the operator in more detail in the next sections. For better understanding the meaning of the $\nu$-operator, see also \cite{GaVV15,GASS16,GaPaSt17,GaPaSt18,GaPaSt21}. 

\subsection{Moschovakis' non-deterministic operator $\nu$}
Let ${\cal A}\in {\sf Struc}$, $Q\subseteq U_{\cal A}^\infty$, and

\vspace{0.2cm}

$\nu[Q](\vec z)=_{\rm df}\{ y_1\in U_{\cal A}\mid ( \vec z \,.\,y_1)\in Q\mbox{ or }(\exists (y_2,\ldots,y_m)\in U_{\cal A}^\infty)(( \vec z \,.\,\vec y) \in Q)\}$. 

\vspace{0.2cm}

\noindent
In the following instructions of type (10), ${\cal O}$ stands for any {\em oracle} ({\em set}) $Q\subseteq U_{\cal A}^\infty$.

\begin{overview}
[Oracle instructions with $\nu$-operator]\label{NuOracleInstr}

\hfill

\nopagebreak 

\noindent \fbox{\parbox{ 11.8cm}{

$\nu$-instructions

\qquad $(10)$ \quad $\ell:\,\, Z_j:=\nu[{\cal O}](Z_1,\ldots,Z_{I_1})$ \hfill (the\, long\, version)

\vspace{0.1cm}

\qquad $(10)$ \quad {\sf $\ell:\,\, Z_j:=\nu[{\cal O}]_{I_1}$} \hfill (the short version)
}}
\end{overview}
Both versions given in Overview \ref{NuOracleInstr} stand for the same instruction. $j$ is again a placeholder for a positive integer. The $\nu$-instructions are strings where the dots $\ldots$ in the long versions of instructions of type (10) are substrings of these instructions. It is also possible to use sets such as $O_Q=\bigcup_{i=1}^\infty \{ (z_1,z_2,\ldots, z_{2i} )\in U_{\cal A}^ {2i}\mid (z_2,z_4,\ldots, z_{2i} ) \in Q\}$ as oracle or, alternatively, we can use pseudo instructions such as the following instructions.

\vspace{0.3cm}

\qquad $^{\rm pseudo}(10)$ \quad $\ell:\,\, Z_j:=\nu[{\cal O}](Z_2,Z_4,\ldots,Z_{I_1})$

\vspace{0.3cm}

\noindent
The execution of a $\nu$-instruction is the process of transforming a configuration by means of a \textcolor{red}{total} multi-valued function presented in Overview \ref{Oper2}. $Q\subseteq U_{\cal A}^\infty$ is fixed. Here, we extend $F_\ell$ and \textcolor{red}{include the case} $\nu[Q] (u_1,\ldots, u_{\nu_1})=\emptyset$ that implies (or means) a looping forever. 

\begin{overview}[Multi-valued maps in ${\cal F}_{\cal M}$ for changing a state]\label{Oper2}
\hfill

\nopagebreak 
\noindent \fbox{\parbox{11.8cm}{
\begin{itemize}\itemsep0pt
\item
$F_\ell: (\bbbn_+\times U_{\cal A} ^{\omega})\,\genfrac{}{}{0pt}{2}{\longrightarrow} {\longrightarrow}\, U_{\cal A} ^{\omega}$ \hfill for instructions of type $(10)$

\parskip 0.5mm 
\begin{tabbing} 
$C_{\ell }(\vec \nu,\bar u)\quad$\=$=( u_1,\ldots,u_{j-1},f_i(u_{j_1},\ldots, u_{j_{m_i}}),u_{j+1}, \ldots)$ \= \kill

$F_\ell=\{((\nu_1,\bar u)\,,\, \, (u_1,\ldots,u_{j-1},y , u_{j+1}, \ldots))\,\,\mid\,\,\, y\in \nu[Q] (u_1,\ldots, u_{\nu_1})\}$\\ 

\hspace{4.23cm}$\textcolor{red}{\cup} \,\, \{((\nu_1,\bar u)\,,\, \bar u)\,\,\mid\,\,\, \nu[Q] (u_1,\ldots, u_{\nu_1})=\emptyset\} $\\

i.e. \,\,\,\,\, $ (\nu_1,\bar u) \genfrac{}{}{0pt}{2}{\longrightarrow} {F_{\ell }} (u_1,\ldots,u_{j-1},y , u_{j+1}, \ldots)$ \,if\,\, $y\in \nu[Q] (u_1,\ldots, u_{\nu_1})$\vspace{0.2cm}
\end{tabbing} 
\end{itemize}}}
\end{overview}

\begin{overview}[Guessing by a $\nu$-oracle machine with operator $\nu$]\label{NuAndGuess}

\hfill

\nopagebreak 

\noindent \fbox{\parbox{11.8cm}{

$ \begin{array}{c}
\hspace{2.05cm}\,z_1\,\, \cdots \,\, \, \,z_{n_0}  \hspace{2.75cm} \,z_1\,\, \cdots \,\, \,\, z_{n_0}\,\,\,\,\,\,\,\,\,\, \, y_1\hspace{1.4cm}\quad\qquad\\
\hspace{1.65cm}\downarrow \qquad\, \,\,\, \,\,\downarrow\, \,\, \hspace{2.8cm}\,\, \downarrow \qquad\,\, \,\,\,\, \downarrow\,\,\, \,\,\, \, \, \,\, \,\,\, \,\, \downarrow\,\, \,\,\, \quad\qquad \hspace{0.3cm}\\
 Z_{j }:= \nu[{\cal O}](Z_1,\ldots,Z_{I_1}); \,\ldots\,; Z_{j}:= \nu[{\cal O}](Z_1,\ldots\textcolor{gray}{,Z_{I_1-1}},Z_{I_1}) \quad \ldots\quad\qquad\\ \!\downarrow\,\,\,\quad\qquad\qquad\quad\qquad\qquad \quad \,\,\,\, \downarrow\qquad\qquad\qquad\qquad\quad\qquad\qquad\quad\quad\qquad\\
\quad\,\,y_1\,\qquad\qquad\qquad\qquad\qquad\quad \,\,\,\, y_2\!\qquad \quad\Rightarrow (z_1,\ldots,z_{n_0},y_1,\ldots,y_m)\in Q
\end{array}$
}}
\end{overview}

\noindent The executions of $m$ $\nu$-instructions allow to guess $m$ values as illustrated in Overview \ref{NuAndGuess}. 

${\sf M}_{\cal A}^\nu(Q)$ is the class of all $\nu$-oracle BSS RAMs over ${\cal A}$ that are able to execute instructions of types $(1)$ to $(8)$ and (10) and to evaluate the oracle $Q\subseteq U_{\cal A}^\infty$ to which $\nu$ may be applied. 

For any ${\cal M}$ in ${\sf M}_{\cal A}^\nu(Q)$, the relations $\to_{\cal M}$ and $(\to_{\cal M})_{{\rm Stop}_{\cal M}}$ are \textcolor{red}{total} or partial multi-valued functions $\to_{\cal M}\,: {\sf S}_{\cal M} \genfrac{}{}{0pt}{2}{\longrightarrow} {\longrightarrow} {\sf S}_{\cal M}$ and $(\to_{\cal M})_{{\rm Stop}_{\cal M}}\!:\,\,\subseteq\! {\sf S}_{\cal M} \genfrac{}{}{0pt}{2}{\longrightarrow} {\longrightarrow} {\sf S}_{\cal M}$ which we therefore also denote by $\genfrac{}{}{0pt}{2}{\longrightarrow} {\longrightarrow}_{\cal M}$ and $(\genfrac{}{}{0pt}{2}{\longrightarrow} {\longrightarrow}_{\cal M}) _{{\rm Stop}_{\cal M}}$. Each configuration assigned by $\to_{\cal M}$ to a configuration $(\ell\,.\,\vec \nu\,.\,\bar u)\in {\sf S}_{\cal M}$ is uniquely defined or defined by applying $ \nu[Q]$ to $ (u_1,\ldots, u_{\nu_1})$ which includes $(\ell\,.\,\vec \nu\,.\,\bar u)\to_{\cal M}(\ell\,.\,\vec \nu\,.\,\bar u)$ for the case that there holds $\nu[Q] (u_1,\ldots, u_{\nu_1})=\emptyset$. ${\rm Input}_{\cal M}$ and ${\rm Output}_{\cal M}$ are defined in the same way as the input and the output procedure for deterministic BSS RAMs. Let ${\rm Input}_{\cal M}(\vec x)= (1\,.\,(n, 1,\ldots, 1)\,.\,(x_1, \ldots,x_n,x_n,x_n, \ldots ))$ for all $\vec x=(x_1,\ldots,x_n)\in U_{\cal A}^\infty$, where $(n, 1,\ldots, 1)\in (\mathbb{N}_+)^{k_{\cal M}}$. Moreover, let ${\rm Output}_{\cal M}(\ell\,.\,\vec \nu\,.\,\bar u)= (u_1,\ldots, u_{\nu_1})$ for any $(\ell\,.\,\vec \nu\,.\,\bar u)\in {\sf S}_{\cal M}$. By (\ref{R3}), ${\rm Res}_{\cal M}(\vec x)$ is the set of all outputs of ${\cal M}$ for the input $\vec x$. For any $Q\subseteq U_{\cal A}^\infty$ and every ${\cal M}$ in ${\sf M}_{\cal A}^\nu(Q)$, the {\em result function} ${\rm Res}_{\cal M}: U_{\cal A}^\infty\to {\mathfrak P}(U_{\cal A}^\infty)$ of ${\cal M}$ is defined by \begin{equation}\tag{R3}\label{R3} \qquad {\rm {\rm Res}_{\cal M}}({\vec x})=\{{\rm Output}_{\cal M}(con)\mid ({\rm Input}_{\cal M}(\vec x),con)\in ( \genfrac{}{}{0pt}{2}{\longrightarrow} {\longrightarrow}_{\cal M})_{{\rm Stop}_{\cal M}}\}
\end{equation} 
for all $\vec x \in U_{\cal A}^\infty$. Let us recall that a partial single-valued function $f:\,\subseteq\! U_{\cal A}^\infty\to U_{\cal A}^\infty$ is non-deterministically {\em $\nu$-computable by a BSS RAM over ${\cal A}$} ({\em by means of the oracle $Q$}) if there is a BSS RAM ${\cal M}$ in ${\sf M}_{\cal A}^\nu(Q)$ such that ${\rm Res}_{\cal M}({\vec x})=\{f(\vec x)\}$ holds if ${\rm Res}_{\cal M}({\vec x})\not=\emptyset$ and $f(\vec x)\uparrow$ holds if $\vec x \in U_{\cal A}^\infty$ and ${\rm Res}_{\cal M}({\vec x})=\emptyset$ hold. A partial multi-valued function $f:\,\subseteq\! U_{\cal A}^\infty \genfrac{}{}{0pt}{2}{\longrightarrow} {\longrightarrow}U_{\cal A}^\infty$ is non-deterministically {\em $\nu$-computable by a BSS RAM over ${\cal A}$} ({\em by means of $Q$}) if there is an ${\cal M}$ in ${\sf M}_{\cal A}^\nu(Q)$ such that $f=\{(\vec x, \vec y)\mid \vec x \in U_{\cal A}^\infty \,\,\&\,\,\vec y\in {\rm Res}_{\cal M}({\vec x})\}$ holds. $ \vec x \genfrac{}{}{0pt}{2}{\longrightarrow} {f} \vec y$ means $(\vec x,\vec y)\in f$. 

For all BSS RAMs over ${\cal A}\in {\sf Struc}$ and all other ${\cal A}$-machines, the terms {\sf result function} and {\sf computable} are defined analogously (cf.\,\,\cite{GASS25A}). A set is \big[ND-\big]{\sf semi-decidable} by a \big[$\nu$-oracle\big] machine if it is the halting set of this machine. A set is {\sf decidable} if it and its complement are semi-decidable by machines in ${\sf M}_{\cal A}$. For more background, see Part I and Part II\,b.

\subsection{Examples}\label{SectionExamples}
The following examples show that it is possible to consider important mathematical functions and operators --- which are sometimes assumed to be computable in a fixed time unit and used, e.g., in computational geometry \cite{PS85} --- and assume their abstract computability in the framework of abstract computation (by BSS RAMs and beyond today's computers). Example \ref{Examp3} shows that the class of $\nu$-computable functions strictly includes the class of computable functions with respect to BSS RAMs over $\mathbb{R}^=$ defined by $\mathbb{R}^==(\mathbb{R};1,0;+,-,\cdot;=)$.

\begin{example}[Computable functions and decidable relations]\label{Example1}\hfill 

\begin{itemize}[leftmargin=0.5cm, itemsep=0.1cm]
\item
The total function $g_1: \mathbb{R}\to \mathbb{R}$ given by $g_1(x)=x^2$ for all $x\in \mathbb{R}$ is computable by an infinite $\mathbb{R}^=$-machine ${\cal M}$ with ${\sf I}_{\cal M}={\sf O}_{\cal M}=\mathbb{R}$, the input procedure given by ${\rm Input}_{\cal M}(x_1)=(1\,.\,{\rm In}_{\cal M}(x_1))=(1\,.\,1\,.\,(x_1,x_1,x_1,\ldots))$, and the program ${\cal P}_{\cal M}$ given by 

\begin{tabular}{l} 
\qquad\qquad {\sf $1: \, Z_1:=f_3^2(Z_1, Z_1);$ \,}\qquad {\sf $2: \,$ stop}.
\end{tabular} \hfill (prog\,1)

\noindent whose execution leads to the output $x_1^2$ by ${\rm Output}_{\cal M}(2\,.\,1\,.\,(x_1^2,x_1,x_1,\ldots)) = x_1^2$. Note, that the partial function $g_2: \subseteq \mathbb{R}^\infty \to \mathbb{R}^\infty$ defined by $g_2(x_1)=x_1^2$ and $g_2(x_1,\ldots,x_n)\uparrow$ for $n>1$ is computable by a BSS RAM over $\mathbb{R}^=$. A BSS RAM in ${\sf M}_{\mathbb{R}^=}$ using also the program given by (prog\,1) computes the total function $g_3: \mathbb{R}^\infty\to \mathbb{R}^\infty$ with $g_3(\vec x)=(x_1^2,x_2, \ldots,x_n)$ for all $\vec x\in \mathbb{R}^\infty$.

\item
The binary relation $Q_1\subseteq \mathbb{R}^\infty$ given by $Q_1=\{(x_1,x_2)\in \mathbb{R}^2\mid x_1=x_2^2\}$ is decidable by a BSS RAM ${\cal M}\in {\sf M}_{\mathbb{R}^=}$ that computes the characteristic function $\chi_{Q_1}: \mathbb{R}^\infty \to \{1,0\}$ --- defined by $\chi_{Q_1}(\vec x)=1$ if $\vec x \in Q_1$ and otherwise $\chi_{Q_1}(\vec x)=0$ (for details see also \cite[pp.\,13--14]{GASS25A} and Part II\,b) --- by using the following input procedure, the program ${\cal P}_{\cal M}$, and the corresponding output procedure used to output $\chi_{Q_1}(\vec x)$ for any input $\vec x\in \mathbb{R}^\infty$. Let ${\rm Input}_{\cal M}(\vec x)=(1\,.\,{\rm In}_{\cal M}(x))=(1\,.\,(n,1)\,.\,(x_1,x_2,\ldots, x_n, x_n,\ldots))$ for $\vec x\in \mathbb{R}^\infty$ and let ${\cal P}_{\cal M}$ be the program

\vspace{0.2cm}

{\sf 
\begin{tabular}{ll}
 \quad $1: \, Z_3:=f_3^2(Z_2, Z_2);$\qquad &\quad\, $5 : \,$ if $I_1=I_2$ then goto $6$ else goto $8;$ \\
\quad $2: \, Z_2:=Z_1;$ &\quad\, $6 : \,$ if $r_1^{2}(Z_2,Z_3)$ then goto $7$ else goto $8;$ \\ 
\quad $3: \, Z_1:=c_2^0;$ &\quad\, $7: \, Z_1:=c_1^0;$ \\ 
\quad $4 : \, I_2=I_2+1;$ &\quad\, $8: \, I_1:=1;$ \\
& \quad\, $9: \,$ stop. \\
\end{tabular}}

\vspace{0.2cm}

\noindent Thus, we get ${\rm Output}_{\cal M}(9\,.\, (1,2)\,.\, (1,x_1,x_2^2, x_2,x_2,\ldots))=1$ if $x_1=x_2^2$ and $\nu_1=2$ hold after creating the initial configuration. Otherwise we get 0. 0 can result from the end configuration $( 9\,.\, (1,2)\,.\, (0,x_1,x_2^2, x_2,x_2, \ldots))$ for any input $(x_1,x_2)\in \mathbb{R}^2$ satisfying $x_1\not=x_2^2$ or from end configurations such as $(9\,.\, (1,2)\,.\, (0,x_1,x_2^2, x_4,x_5, x_5,\ldots))$ for an input of length $5$, and so on. 

\item The set $Q_2$ given by $Q_2=\{(x_1,x_2, x_3)\in \mathbb{R}^3\mid x_1=x_2^2\,\,\&\,\, x_2=x_3^2\}$ is decidable by a BSS RAM ${\cal M}\in {\sf M}_{\mathbb{R}^=}$ by means of the following program ${\cal P}_{\cal M}$. 

\vspace{0.2cm}

{\sf 
\begin{tabular}{ll}
\quad $1: \, Z_4:=f_3^2(Z_3, Z_3);$\quad&\quad\, $7 : \,$ if $I_1=I_2$ then goto $8$ else goto $11;$ \\
\quad $2: \, Z_3:=f_3^2(Z_2, Z_2);$ &\quad\, $8 : \,$ if $r_1^2(Z_5,Z_3)$ then goto $9$ else goto $11;$ \\
\quad $3: \, Z_5:=Z_1;$ &\quad\, $9 : \,$ if $r_1^2(Z_2,Z_4)$ then goto $10$ else goto $11;$ \\
\quad $4: \, Z_1:=c_2^0;$ &\quad $10: \, Z_1:=c_1^0;$\\ 
\quad $5 : \, I_2=I_2+1;$ &\quad $11: \, I_1:=1;$ \\
\quad $6 : \, I_2=I_2+1;$ &\quad $12: \,$ stop.\\
\end{tabular}}

\vspace{0.2cm}

\end{itemize}
\end{example}

\begin{example}[Non-deterministically computable inverses]\hfill 

\begin{itemize}[leftmargin=0.5cm, itemsep=0.1cm]
\item Let $g_4$ be the partial single-valued function $g_4:\,\,\subseteq \! \mathbb{R}^\infty \to \mathbb{R}^\infty$ given by $g_4(\vec x)=\sqrt{x_1}$ for all $\vec x \in \mathbb{R}^1 $ with $x_1\geq 0$ and otherwise $g_4(\vec x)\uparrow$. This means that $g_4(\vec x)$ is not defined for $\vec x \in \mathbb{R}^1 $ with $x_1< 0$ and for $\vec x \in \mathbb{R}^\infty \setminus \mathbb{R}^1$. $g_4$ is {\rm ND}-computable by a machine ${\cal M}\in {\sf M}^{\rm ND}_{\mathbb{R}^=}$ whose input procedure provides initial configurations by $\vec x \,\genfrac{}{}{0pt}{2}{\longrightarrow} {\,{\rm Input }_{\cal M}} (1\,.\,(n,1)\,.\,(x_1,\ldots, x_n, y_1,\ldots,y_m, x_n,x_n, \ldots))$ --- such that $(\vec x, (1,n,1,x_1,\ldots, x_n, y_1,\ldots,y_m, x_n,x_n, \ldots))\in {\rm Input}_{\cal M}$ holds for any $m\geq 1$ and any guesses in $(y_1,\ldots,y_m)\in \mathbb{R}^\infty$ ---
and whose program ${\cal P}_{\cal M}$ is given by
\vspace{0.2cm}

\begin{tabular}{l}
\qquad\qquad {\sf $1 : \,$ if $I_1=I_2$ then goto $2$ else goto $1;$ }\\

\qquad\qquad {\sf $2: \, Z_3:=f_3^2(Z_2, Z_2);$ }\\

\qquad\qquad {\sf $3 : \,$ if $r_1^2(Z_1,Z_3)$ then goto $4$ else goto $3;$ }\\

\qquad\qquad {\sf $4: \, Z_3:=f_3^2(Z_4,Z_4);$ }\\

\qquad\qquad {\sf $5 : \,$ if $r_1^2(Z_2,Z_3)$ then goto $6$ else goto $5;$ }\\

\qquad\qquad {\sf $6: \, Z_1:=Z_2;$ }\\

\qquad\qquad {\sf $7: \,$ stop}.
\end{tabular}

\vspace{0.2cm}

\noindent If $c(I_1)=1$, $c(Z_1)=x_1 \geq 0$, $m\geq 4$, $c(Z_2)=y_1=\sqrt{x_1}$, and $c(Z_4)=y_3=\sqrt[4]{x_1}$ hold after creating an initial configuration, then the end configuration $con$ is $(7\,.\,(1,1)\,.\, (\sqrt{x_1},\sqrt{x_1}, \sqrt{x_1}, \sqrt[4]{x_1}, y_4,\ldots, y_m, x_1,x_1,\ldots))$ which implies the output ${\rm Output}_{\cal M}(con) = \sqrt{x_1}$. The result function defined by ${\rm Res_{\cal M}}(\vec x)=\{ {\rm Ouput}_{\cal M}( (\to_{\cal M})_{{\rm Stop}_{\cal M}} (1\,.\,\vec \nu\,.\,\bar u))\mid (\vec x, (\vec \nu\,.\,\bar u))\in {\rm In}_{\cal M}\}$ is given by 

\vspace{0.1cm}

\qquad ${\rm Res}_{\cal M}(\vec x)= \left\{\begin{array}{ll} \{\sqrt{x_1}\} & \mbox{ if } \vec x \in \mathbb{R}^1 \,\,\&\,\, x_1\geq 0,\\

\emptyset & \mbox{ otherwise}.\end{array}\right.$ \hfill (res\,1)

\item
The partial multi-valued function $g_5:\,\subseteq \! \mathbb{R}^\infty \genfrac{}{}{0pt}{2}{\longrightarrow} {\longrightarrow} \mathbb{R}^\infty$ given by $g_5=\{(x, y)\mid x \in \mathbb{R} \,\,\&\,\,y\in \{\sqrt{x},-\sqrt{x}\}\}$ is non-determin\-istically computable by a machine ${\cal M}\in {\sf M}^{\rm ND}_{\mathbb{R}^=}$ using the following program ${\cal P}_{\cal M}$.

\vspace{0.2cm}

\begin{tabular}{l}

\qquad\qquad {\sf $1 : \,$ if $I_1=I_2$ then goto $2$ else goto $1;$ }\\

\qquad\qquad {\sf $2: \, Z_3:=f_3^2(Z_2,Z_2);$ }\\

\qquad\qquad {\sf $3 : \,$ if $r_1^2(Z_1,Z_3)$ then goto $4$ else goto $3;$ }\\

\qquad\qquad {\sf $4: \, Z_1:=Z_2;$ }\\

\qquad\qquad {\sf $5: \,$ stop}.
\end{tabular}

\vspace{0.2cm}

$(1,1,1,x_1,\sqrt{x_1}, y_2, x_1, \ldots)\to_{\cal M} ^*(5,1,1,\sqrt{x_1}, \sqrt{x_1}, x_1,x_1, \ldots)$ holds for the transitive closure $\to_{\cal M} ^*$ of the relation $\to_{\cal M} $ and $m=2$. Moreover, we have $(1,1,1,x_1,-\sqrt{x_1}, y_2, y_3,x_1, \ldots)\to_{\cal M} ^*(5,1,1,-\sqrt{x_1}, -\sqrt{x_1}, x_1, y_3,x_1,x_1, \ldots)$ for $m=3$, and so on. Consequently, we obtain (res\,2) which implies $\vec x \genfrac{}{}{0pt}{2}{\longrightarrow} {g_5} \sqrt{x_1}$ and $\vec x \genfrac{}{}{0pt}{2}{\longrightarrow} {g_5}- \sqrt{x_1}$ for all $\vec x \in \mathbb{R}^1 $ with $x_1\geq 0$. 

\vspace{0.1cm}

\qquad ${\rm Res}_{\cal M}(\vec x)= \left\{\begin{array}{ll} \{\sqrt{x_1},-\sqrt{x_1}\} & \mbox{ if } \vec x \in \mathbb{R}^1 \,\,\&\,\, x_1\geq 0,\\

\emptyset & \mbox{ otherwise}\end{array}\right.$ \hfill (res\,2)

\vspace{0.1cm}

\end{itemize}
\end{example}

\begin{example}[Non-deterministically $\nu$-computable inverses]\label{Examp3}\hfill 

\noindent For the decidable set $Q_2$ given by $Q_2=\{(x_1,y_1, y_2)\in \mathbb{R}^3\mid x_1=y_1^2\,\,\&\,\, y_1=y_2^2\}$, we have $\nu[Q_2](x_1)=\{\sqrt{x_1}\}$ for all $x_1\geq 0$. $\nu[Q_2](x_1,x_2)=\{\sqrt{x_2}, -\sqrt{x_2}\}$ holds for all $(x_1,x_2)$ satisfying $x_2\geq 0$ and $x_2^2=x_1$. For $\vec x \in \mathbb{R}^3$, we have $\nu[Q_2](\vec x)=\emptyset$. 

\begin{itemize}[leftmargin=0.5cm, itemsep=0.1cm]
\item The result function ${\rm Res}_{\cal M}$ of the $\nu$-oracle machine ${\cal M}\in {\sf M}_{\mathbb{R}^=}^\nu(Q_2)$ using its program ${\cal P}_{\cal M}$ given by 

\vspace{0.2cm}

\begin{tabular}{l}

\qquad\qquad {\sf $1 : \,$ if $I_1=I_2$ then goto $2$ else goto $1;$}\\

\qquad\qquad {\sf $2: \, Z_1:= \nu[{\cal O}](Z_1,\ldots,Z_{I_1});$}\\

\qquad\qquad {\sf $3: \,$ stop}.
\end{tabular}\hfill (prog\,2)

\vspace{0.2cm}

\noindent is equal to the result function given by (res\,1). Hence, $g_4$ is $\nu$-computable.

\item Any function $g$ satisfying $g(\vec x)=\sqrt{x_1}$ for all $\vec x \in \mathbb{R}^1$ with $x_1\geq 0$ is not computable by a BSS RAM in ${\sf M}_{\mathbb{R}^=}$ since for any ${\cal N}\in {\sf M}_{\mathbb{R}^=}$ and any $n\in \bbbn$, there holds ${\rm Res}_{\cal N}(n)\in \bbbn^\infty$. 

\item By deleting the first instruction in (prog\,2), we can obtain the program ${\cal P}_{{\cal M}'}$ of a machine ${\cal M}'\in {\sf M}_{\mathbb{R}^=}^\nu(Q_2)$ whose  result function is given by (res\,3).

\vspace{0.1cm}

\qquad ${\rm Res}_{{\cal M}'}(\vec x)= \left\{\begin{array}{ll} \{\sqrt{x_1}\} & \mbox{ if } \vec x \in \mathbb{R}^1 \,\,\&\,\, x_1\geq 0,\\

\{\sqrt{x_2},-\sqrt{x_2}\} & \mbox{ if } \vec x\in \mathbb{R}^2 \mbox{ and } x_2^2=x_1,\\

\emptyset & \mbox{ otherwise.}\end{array}\right.$ \hfill (res\,3)

\vspace{0.1cm}
\end{itemize}
\end{example}

\subsection{The power of non-deterministic BSS RAMs}

Let $({\rm SDEC}_{\cal A}^\nu)^Q$ be the class of all problems non-deterministically semi-decidable by $\nu$-oracle machines in ${\sf M}_{\cal A}^\nu(Q)$. The ideas leading to Proposition 10 in \cite[p.\,592]{GASS20} (and Proposition 1 in \cite{GASS19Preprint}, respectively) can also be used to prove the following statements which include results presented in 2015 (for details see \cite{GaVV15}). We think that this helps to understand the similarity between the hierarchy introduced by Moschovakis in \cite{MOSCHO} and the hierarchy that will be considered in Part IV. For the definition of 3-tape machines, see Overview \ref{NDThreeTape}. For more details and the terms {\sf subprogram} and {\sf pseudo instruction}, see Section \ref{SectSimMachinesWithNuInst}. 

\begin{theorem}[$\nu$-queries and the ${\rm ND}$-semi-decidability]\label{OracleMNondetM} Let ${\cal A}\in {\sf Struc}$, $Q\subseteq U_{\cal A}^\infty$, and $Q\in {\rm SDEC}_{\cal A}$. Then, there holds \[({\rm SDEC}_{\cal A}^\nu)^Q\subseteq {\rm SDEC}_{\cal A}^{\rm ND}.\]
 \end{theorem}
{\bf Proof.} 
Roughly speaking, a BSS RAM ${\cal M}\in {\sf M}_{\cal A}^{\nu}(Q)$ can be simulated by a non-deterministic multi-tape ${\cal A}$-machine ${\cal N}_{\cal M}\in ({\sf M}_{\cal A}^{(3)})^{\rm ND}$ --- including the stop instruction of type (8) --- as follows. 

Let ${\cal N}_Q\in {\sf M}_{\cal A}$ be one of the BSS RAMs in $ {\sf M}_{\cal A}$ that semi-decide $Q$. Each execution of a $\nu$-instruction in ${\cal P}_{\cal M}$ can be simulated by using a subprogram derived from the program ${\cal P}_{{\cal N}_Q}$. For this purpose, let ${\cal N}_{\cal M}$ be equipped with three tapes: a first tape given by $Z_{1,1},Z_{1,2},\ldots$ for storing all individuals of a start configuration of ${\cal N}_{\cal M}$ --- including the input $(x_1,\ldots,x_n)$ for ${\cal N}_{\cal M}$ and guesses --- during the entire runtime of the simulation, a second tape $Z_{2,1},Z_{2,2},\ldots$ for storing the values $c(Z_1),c(Z_2),\ldots$ of the $Z$-registers of ${\cal M}$ and for simulating the execution of all instructions of types (1) to (8) in ${\cal P}_{\cal M}$ by executing instructions of a subprogram ${\cal P}_{\cal M}'$ --- derived from ${\cal P}_{\cal M}$ --- in the analogous way as ${\cal M}$ on the input given for ${\cal N}_{\cal M}$, and a third tape $Z_{3,1},Z_{3,2},\ldots$ for evaluating all $\nu$-instructions, executed by ${\cal M}$, by using a suitable copy ${\cal P}''_{{\cal N}_Q}$ of ${\cal P}_{{\cal N}_Q}$. We assume that each $n$-ary input for ${\cal N}_{\cal M}$ is also the input for $\cal M$ and stored in $Z_{1,1},Z_{1,2},\ldots, Z_{1,n}$. The guesses $y_1,\ldots, y_{m_1}$, $y_{m_1+1},\ldots, y_{m_1+m_2}$, $y_{m_1+m_2+1},\ldots, y_{m_1+m_2+m_3},\ldots $ of ${\cal N}_{\cal M}$ are then stored in its registers $Z_{1,n+1}, Z_{1,n+2},\ldots$. Let $m_0$ and $i_0$ be $0$ at the beginning of the simulation of ${\cal M}$. If, in simulating ${\cal M}$, 

\begin{itemize}\parskip -0.5mm 
\item the $(i_0+1)^{\rm th}$ $\nu$-instruction --- in the sequence of instructions performed by ${\cal M}$ during the computation process --- is under consideration, 
\item the content $c(I_1)$ of the register $I_1$ of ${\cal M}$ --- which is then also stored in $I_{2,1}$ and $I_{3,1}$ --- is $n_0$, and 
\item $(z_1,\ldots, z_{n_0})$ is the input $(c(Z_1),\ldots, c(Z_{n_0}))$ of the $(i_0+1)^{\rm th}$ $\nu$-instruction,
\end{itemize}
then let $m_0=\sum_{i=1}^{i_0}m_i$ (where \textcolor{blue}{$\sum_{i=1}^{0}m_i=_{\rm df}0$}) and let ${\cal N}_{\cal M}$ evaluate this instruction (without simulating the stop instruction of ${\cal N}_Q$) as follows. 

\vspace{0.1cm}

{\sf 
\hspace{0.5cm}for $s:=1,2,\ldots$ do \{ $m:=ca_1(s); $ $s_0:=ca_2(s); $ 

\hspace{1.1cm}if $m=0$ or $s_0=0$ then $m:=1;$ $s_0:=1;$

\hspace{1.15cm}simulate $s_0$ steps of ${\cal N}_Q$ on input $(z_1,\ldots, z_{n_0},y_{m_0+1},\ldots, y_{m_0+m})$

\hspace{0.92cm}\begin{tabular}{l}and \textcolor{blue}{ if the label $\ell_{{\cal P}_{{\cal N}_Q}}$ is reached} then \\\hspace{0.55cm}  execute \{$z:=y_{m_0+1};$ $m_0:= m_0+m;$ output $z$\} 
\quad \end{tabular} \hspace{0.3cm}\textcolor{blue}{{\it (ps.instr\,1)}}

\hspace{1.12cm}\}} \hfill {\it (algo\,1)}

\vspace{0.1cm}

\noindent 
Here, we use a loop construct allowing to repeat the simulation in the block between the outer braces for $s:=1,2,\ldots$ (and for $m$ and $s_0$ that generally result from $s=\frac12((m+s_0)^2 +3s_0+m)$). The instructions in the block allow to simulate --- at most --- $\max\{1,ca_2(s)\}$ steps of ${\cal N}_Q$. \textcolor{blue}{\it (ps.instr\,1)} means that the execution of the block is repeated for $s:=1,2,\ldots$ {\sf only until} the label $\ell_{{\cal P}_{{\cal N}_Q}}$ of the stop instruction of ${\cal N}_Q$ is reached for the first time in the process of simulating ${\cal N}_Q$. Moreover, this also means that the execution of this loop construct could be repeated forever and thus that, for the used sequence of guesses $y_{m_0+1}, y_{m_0+2},\ldots$, ${\cal N}_{\cal M}$ {\sf loops forever}. After recognizing the label $\ell_{{\cal P}_{{\cal N}_Q}}$, the subprogram {\sf output $z$} can be executed. It causes that the loop construct is exited and the guess $y_{m_0+1}$ is returned. If ${\cal N}_Q$ halts on $(z_1,\ldots, z_{n_0},y_{m_0+1},\ldots, y_{m_0+m}) $, then $y_{m_0+1}$ is one of the possible output values that can be the result of the execution of the $\nu$-instruction --- i.e. one of the theoretically conceivable values --- and that can be assigned to the $Z$-register $Z_j$ by ${\cal M}$ and that is therefore assigned to the $Z$-register $Z_{2,j}$ of ${\cal N}_{\cal M}$. In this case, let $m_{i_0+1}$ be the value $m$ that is determined in the loop (for the last considered $s$) in which the simulation comes to evaluate the label $\ell_{{\cal P}_{{\cal N}_Q}}$. After this, i.e., after simulating the execution of a $\nu$-instruction {\sf instruction$_k$} in ${\cal P}_{\cal M}$, ${\cal N}_{\cal M}$ continues the simulation of ${\cal M}$ by evaluating the instruction {\sf instruction$_{k+1}$} in ${\cal P}_{\cal M}$. For this reason, let us also carry out $i_0:=i_0+1$ so that we can again use the new $i_0$ in the description of a further $\nu$-instruction. Consequently, we can say that the guesses $y_{\sum_{i=1}^{i_0} m_i+1}$, \ldots, $y_{\sum_{i=0}^{i_0+1} m_i} $ are used later in simulating the execution of the $(i_0+1)^{\rm th}$ $\nu$-instruction (for the incremented $i_0$). If ${\cal N}_{\cal M}$ would need more than its own guesses in the attempt to simulate ${\cal M}$, then ${\cal N}_{\cal M}$ uses also the last component $x_n$ of its own input in the following. If ${\cal M}$ executes the stop instruction after $t$ steps, then such a computation can be simulated by using $\sum_{i=0}^t m_i$ suitable guesses of ${\cal N}_{\cal M}$. Analogously to the proof of Proposition 10 in \cite{GASS20}, we can show that there is a machine in ${\sf M}_{\cal A}^{\rm ND}$ that can simulate the multi-tape machine ${\cal N}_{\cal M}$. 
\qed

\begin{remark}[Another algorithm with identity test]\label{RemOracleMNondetM}

The $\nu$-instructions \linebreak can help to find answers to the following questions for $\vec z=(z_1,\ldots, z_{n_0})$. Is there a tuple of the form $(z_1,\ldots, z_{n_0}, y_1,\ldots,y_k)$ (for $k\geq1$) in a given oracle set $Q$? Which components can follow the first components $z_1,\ldots,z_{n_0}$ in a tuple that belongs to $Q$? 

The process of producing guesses by $\nu$-instructions for answering the second question is only useful if ${\cal A}$ contains at least two elements. Let ${\cal A}$ contain at least two elements and the identity relation over $U_{\cal A}$. Then, instead of extending the sequence of $m$ guesses gradually as described by {\rm {\it (algo\,1)}}, it is also possible to use only one of the following sequences of guesses,

\vspace{0.1cm}

\hspace{1.4cm} $y_1,y_3,\ldots, y_{2m_1-1}$ 

\qquad or \quad $y_{2m_1+1},y_{2m_1+3}, \ldots, y_{2m_1+2m_2-1}$ 

\qquad or \quad $y_{2m_1+2m_2+1}$, $y_{2m_1+2m_2+3},\ldots, y_{2m_1+2m_2+2m_3-1} $

 \qquad or \quad $\ldots$,

\noindent as additional input components for ${\cal N}_Q$. In evaluating one of the $\nu$-instructions as described below, the length of the used sequence results from 

\vspace{0.1cm}

\hspace{1.4cm} $z_{n_0} =y_2=y_4=\cdots = y_{2m_1-2}\not= y_{2m_1}$, \footnote{That means that $z_{n_0}\not= y_{2m_1}$ if $m_1=1$, and so on.}

 \qquad or \quad  $ y_{2m_1}=y_{2m_1+2}=y_{2m_1+4}=\cdots = y_{2m_1+2m_2-2}\not= y_{2m_1+2m_2}$, 
\ldots .

\vspace{0.2cm}

\noindent Let $y_0=z_{n_0}$. The evaluation of the $(i_0+1)^{\rm th}$ $\nu$-instruction is possible as follows. 

\vspace{0.1cm}

{\sf 
\hspace{1.35cm} $m:=1;$

\hspace{0.9cm} $\ell:$ if $y_{2m_0+2m}\not=y_{2m_0+2m-2}$ then goto $\ell+2$ else goto $\ell+1;$

\hspace{0.3cm} $\ell+1:$ $m:=m+1;$ goto $\ell;$ 

\hspace{0.3cm} $\ell+2:$ simulate ${\cal N}_Q$ on input $(z_1,\ldots, z_{n_0},y_{2m_0+1},\ldots, y_{2m_0+2m-1}) $ 

\hspace{1.3cm} only until the label $\ell_{{\cal P}_{{\cal N}_Q}}$ is reached$;$ 

\hspace{1.3cm} $z:=y_{2m_0+1};$ $m_0:= m_0+m;$ output $z$} \hfill {\rm {\it (algo\,2)}}
\end{remark}
\begin{theorem}[Guessing and the $\nu$-semi-decidability]\label{OracleMEqNondM} Let ${\cal A}\in {\sf Struc}$ and $Q=U_{\cal A}^\infty$. Then, there holds \[({\rm SDEC}_{\cal A}^\nu)^Q={\rm SDEC}_{\cal A}^{\rm ND}.\] 
\end{theorem}
{\bf Proof.} Because of $U_{\cal A}^\infty\in {\rm SDEC}_{\cal A}$ and Theorem \ref{OracleMNondetM},   it is enough to show the inclusion ``$\supseteq$''  and that the running of ${\cal M}\in {\sf M}_{\cal A}^{\rm ND}$ on  $\vec x$ can be simulated by a 1-tape $\nu$-oracle machine ${\cal N}_{\cal M}$ that uses its registers $Z_1,Z_2,\ldots$ and $I_1,\ldots, I_{k_{\cal M}}$ for storing the values of the registers of ${\cal M}$ and further registers such as $ I_{k_{\cal M}+1}$ whose initial value $c(I_{k_{\cal M}+1})$ should be $n$. By executing the subprogram

\vspace{0.05cm}

\hspace{1.35cm} $I_{k_{\cal M}+1}:=I_{k_{\cal M}+1}+1 $; \, $Z_{I_{k_{\cal M}+1}}:=\nu[{\cal O}](Z_{1},\ldots, Z_{I_{1}})$, 

\vspace{0.1cm}

\noindent ${\cal N}_{\cal M}$ can generate guesses before it uses them like the guesses of ${\cal M}$. By applying the $\nu$-operator, ${\cal N}_{\cal M}$ receives a value $y_i$ that can be copied into one of the registers $Z_{n+1}, Z_{n+2},\ldots$, before ${\cal M}$ would use the corresponding register for the first time. Note, that the guessed content $c(Z_{n+i})$ can also be equal to $x_n$. To assign $y_i$ to one of the mentioned registers, $I_{k_{\cal M}+1}$ is used to count the registers that could be used in the worst case by ${\cal M}$ until that point. The process of copying $y_i$ into a suitable register $Z_{c(I_{k_{\cal M}+1})}$ can be done by a subprogram ${\sf initguess}(Z_{I_{k_{\cal M}+1}})$ listed in Overview \ref{PseudoInst} and is to be carried out only if an instruction of type (7) is to be simulated. If ${\cal N}_{\cal M}$ would halt on an input after the simulation of $t$ steps performed by ${\cal M}$, then at most $t$ guesses need to be generated by ${\cal N}_{\cal M}$. If ${\cal N}_{\cal M}$ would try to copy an infinite number of guesses into further $Z$-registers in the attempt to simulate ${\cal M}$, then ${\cal M}$ and ${\cal N}_{\cal M}$ do not halt. 
\qed

\vspace{0.2cm}

Since the proofs of Theorems \ref{OracleMNondetM} and \ref{OracleMEqNondM} are based on the possibility to simulate mutually the considered machines --- instruction by instruction --- and to provide or use any combination of guesses, the statements can be generalized. However, the output must be copied onto the first tape (and so on, if necessary) after a simulation and  before  the output procedure can be performed.

\begin{theorem}[Non-deterministic computability]\label{NonDetClasses1}
Let ${\cal A}$ be in $ {\sf Struc}$, $Q\subseteq U_{\cal A}^\infty$, and $Q\in {\rm SDEC}_{\cal A}$. Then, $\{{\rm Res}_{\cal M}\mid {\cal M}\in {\sf M}_{\cal A}^\nu(Q)\}\subseteq \{{\rm Res}_{\cal M}\mid {\cal M}\in {\sf M}_{\cal A}^{\rm ND}\}$ holds. For  $Q=U_{\cal A}^\infty$, there holds $\{{\rm Res}_{\cal M}\mid {\cal M}\in {\sf M}_{\cal A}^\nu(Q)\}= \{{\rm Res}_{\cal M}\mid {\cal M}\in {\sf M}_{\cal A}^{\rm ND}\}$.
\end{theorem}

\subsection{Evaluation of $\nu$-instructions for $n$-ary or decidable $Q$'s}
For the following special cases, we want to give simpler proofs.
\begin{corollary}[Special oracles]\label{OracleMNondetM1} 
Let ${\cal A}$ be a structure with identity relation. 
\[({\rm SDEC}_{\cal A}^\nu)^Q\subseteq {\rm SDEC}_{\cal A}^{\rm ND}\]
holds, if we have {\rm (a1)} or {\rm (a2)} or {\rm (a3)}. 
\begin{enumerate}[label={\rm \hspace*{1cm} (a\arabic*)}]\parskip -0.5mm 
\item $Q\subseteq U_{\cal A}^{n_Q}$ for some $n_Q\geq 1$\hspace{0.3cm}$\&$\hspace{0.3cm}$Q\in\hspace{0.2cm}{\rm DEC}_{\cal A}$\hspace{0.3cm}$\&$\hspace{0.3cm}${\cal A}\in{\sf Struc}_{c_1,c_2}$,
\item $Q\subseteq U_{\cal A}^{n_Q}$ for some $n_Q\geq 1$\hspace{0.3cm}$\&$\hspace{0.3cm}$Q\in {\rm SDEC}_{\cal A}$\hspace{0.3cm}$\&$\hspace{0.3cm}${\cal A}\in{\sf Struc}$,
\item $Q\subseteq U_{\cal A} ^\infty$\hspace{3.1cm}$\&$\hspace{0.3cm}$Q\in\hspace{0.17cm}{\rm DEC}_{\cal A}$\hspace{0.3cm}$\&$\hspace{0.3cm}${\cal A}\in{\sf Struc}_{c_1,c_2}$.
\end{enumerate} 
\end{corollary}
{\bf Proof.} In cases (a1) and (a3), let ${\cal N}_Q\in {\sf M}_{\cal A}$ be any suitable machine that decides the corresponding oracle set $Q$ by computing $\chi_Q: U_{\cal A}^\infty\to \{c_1,c_2\}$ (cf.\, Example \ref{Example1}). For more details, see Part II\,b. In case (a2), let ${\cal N}_Q$ be any machine in ${\sf M}_{\cal A}$ semi-deciding $Q$. In each of these cases, let ${\cal M}$ be any machine in ${\sf M}_{\cal A}^{\nu}(Q)$. Then, we can define a non-deterministic 3-tape ${\cal A}$-machine ${\cal N}_{\cal M}$ for simulating ${\cal M}$. In all cases, let ${\cal N}_{\cal M}$ be again equipped with the three tapes as described above. ${\cal N}_{\cal M}$ uses $Z_{3,1},Z_{3,2},\ldots$ for evaluating the $\nu$-instructions of ${\cal M}$ by using the suitable program ${\cal P}_{{\cal N}_Q}$. The input for ${\cal N}_{\cal M}$ is the input for ${\cal M}$ and stored in $Z_{1,1},\ldots, Z_{1,n}$. The guesses $y_1,\ldots, y_{m_1}$, $y_{m_1+1},\ldots, y_{m_1+m_2}$, $y_{m_1+m_2+1},\ldots, y_{m_1+m_2+m_3},\ldots $ of ${\cal N}_{\cal M}$ are stored in its registers $Z_{1,n+1}$, $Z_{1,n+2}$, $\ldots$. The integers $m_0$ and $i_0$ are equal to $0$ at the beginning of the simulation of ${\cal M}$. If $(z_1,\ldots, z_{n_0})$ is the input $(c(Z_1),\ldots, c(Z_{c(I_1)}))$ of the $(i_0+1)^{\rm th}$ $\nu$-instruction under consideration, then ${\cal N}_{\cal M}$ simulates the execution of this instruction as described below for the three cases individually. Each of the following constructs can only be exited (to evaluate the next instruction in ${\cal P}_{\cal M}$) when the subprogram {\sf output $z$} can be executed after reaching the label $\ell_{{\cal P}_{{\cal N}_Q}}$ of the stop instruction of ${\cal N}_Q$. In the cases (a1) and (a3),  the pseudo instruction  {\sf if $Z_j=c_1$ then goto $\ell_1$ else goto $\ell_2$} allows to decide $\{c_1\}$.

\vspace{0.3cm}

\begin{enumerate}[label={(a\arabic*)}]
\item\parskip -0.2mm 
{\sf 

 If $n_0< n_Q$ then \{ 

\hspace{2.5cm} $m:= n_Q-n_0;$

\hspace{2.5cm} simulate ${\cal N}_Q$ on input $(z_1,\ldots, z_{n_0},y_{m_0+1},\ldots, y_{m_0+m}); $

\hspace{2.5cm} if the output of ${\cal N}_Q$ is $c_1$

\hspace{2.9cm} then \{ $z=y_{m_0+1};$ $m_0:= m_0+m;$ output $z$ \}

\hspace{2.9cm} {\sf else goto $\ell;$}

\hspace{2.4cm} \}}

\hspace{1.6cm} {\sf else goto $\ell;$}

 $\ell:\,$ {\sf goto $\ell$}

\vspace{0.3cm}

\item {\sf If $n_0< n_Q$ then \{ 

\hspace{2.5cm} $m:= n_Q-n_0;$ 

\hspace{2.5cm} simulate ${\cal N}_Q$ on input $(z_1,\ldots, z_{n_0},y_{m_0+1},\ldots, y_{m_0+m}); $ 

\hspace{2.5cm} $z=y_{m_0+1};$

\hspace{2.5cm} $m_0:= m_0+m;$ 

\hspace{2.5cm} output $z$; 

\hspace{2.4cm} \}}

\hspace{1.6cm} {\sf else goto $\ell;$} 

 $\ell:\,$ {\sf goto $\ell$}   \hfill {\scriptsize (Note, the identity relation is not necessary in this case.)}

\item {\sf for $m:=1,2,\ldots$ do 

\hspace{2.4cm} \{

\hspace{2.5cm} simulate ${\cal N}_Q$ on input $(z_1,\ldots, z_{n_0},y_{m_0+1},\ldots, y_{m_0+m}); $

\hspace{2.5cm} if the output of ${\cal N}_Q$ is $c_1$ 

\hspace{2.9cm} then \{ $z=y_{m_0+1};$ $m_0:= m_0+m;$ output $z$ \}
 
\hspace{2.4cm} \}} 
\end{enumerate}
\noindent If the considered machine ${\cal M}$ executes its stop instruction after $t$ steps, then such a computation can be simulated by using $\sum_{i=0}^t m_i$ guesses of ${\cal N}_{\cal M}$. 
\qed

\section{Simulation of $\nu$-oracle machines}\label{SectSimMachinesWithNuInst}

For a detailed description of the proof of Theorem \ref{OracleMNondetM}, we introduce some subprograms. These and similar subprograms are used in the following.

\subsection{Useful subprograms and more}

Let  ${\cal L}$ be a set $ \{\ell_1,\ldots,\ell_k\}$  whose $k$ elements we call {\sf labels} and let {\sf subpr} be a string of the form

\vspace{0.15cm}

${\sf instruction}_1'; \quad \ell_2: {\sf instruction}_2'; \quad \ldots\,;\quad \ell_{k-1}: {\sf instruction}_{k-1}'$ 

\vspace{0.15cm}

\noindent derived from a program ${\cal P}$ of the form \textcolor{blue}{(*)} (see p.\,\,\pageref{star}) by replacing each label $i\in \{1,\ldots,k\}$ (we assume that $k=\ell_{\cal P}$ holds) by $\ell_i\in {\cal L}$ such that, for all $j\in \{1,\ldots,k-1\}$, every instruction {\sf instruction}$_j'$ results from {\sf instruction}$_j$ in ${\cal P}$ by replacing any label $i$ by  $\ell_i$. For simplicity, we call the string of the form 

\vspace{0.15cm}

$\ell_1: {\sf subpr}; \, \ell_k: {\sf stop.}$ 

\vspace{0.15cm}

\noindent  also a {\em program} --- in a broader sense ---  {\em with labels in ${\cal L}$}. Any string 
\vspace{0.15cm}

$\ell_i: {\sf instruction}_{i}'; \quad \ldots\,;\quad \ell_j: {\sf instruction}_j';\,\,\,\,\big[\ell_{j+1}: {\sf stop.}\big]$ \hfill{\small($1\leq i\leq j<k$)}

\vspace{0.15cm}

\noindent   is
a {\em program segment}. 
{\sf subpr} can be a ({\em proper}) {\em  subprogram} of a program ${\cal P}$ if 

\vspace{0.15cm}

$\ell: {\sf subpr};$ 

\vspace{0.15cm}

\noindent is a program segment of ${\cal P}$  (which means that any label of an instruction in {\sf subpr} is not the label of an instruction in another part of ${\cal P}$). Here, ${\cal P}$  can be a program in a broader sense. For further possible definitions see \cite{GASS20}. 

We generally use many subprograms that can be summarized and described by so-called pseudo instructions. The {\em pseudo instructions} are also something like the names for subprograms or for other (sequences of) instructions. Some of them are given in Overview \ref{PseudoInst} and in Overviews 2, 4, and 5 in \cite{GASS20}.

\begin{overview}
[Some pseudo instructions]\label{PseudoInst}

\hfill

\nopagebreak 

\noindent \fbox{\parbox{11.8cm}{

\quad {\sf goto} $\ell$ 

\vspace{0.1cm}

\quad $(Z_{2,1},\ldots, Z_{2,I_{2,1}}):= (Z_{1,1},\ldots, Z_{1,I_{1,1}})$

\vspace{0.1cm}

\quad {\sf if $I_{1,2}=i$ then goto $i\,'$}

\vspace{0.1cm}

\quad {\sf for $I_{1,3}:=1,2,\ldots \big[,I_{1,4}$\big] do inst}

\vspace{0.1cm}

\quad $ I_{1,4}:=ca_1^{\big[+\big]}(I_{1,3});$ $I_{1,5}:=ca_2^{\big[+\big]}(I_{1,3})$

\vspace{0.1cm}

\quad $I_{1,1}:=I_{1,1} +I_{1,4}$

\vspace{0.1cm}

\quad ${\sf init}(Z_{I_{k_{\cal M}+1}})$ \, and \, ${\sf init}(Z_{d_1,I_{d_1,k_{\cal M}+1}})$ \, and \, 
${\sf initguess}(Z_{I_{k_{\cal M}+1}})$ 

\vspace{0.1cm}

\quad {\sf if $I_{1,9}=I_{1,5}$ then inst$_1$ \big[else inst$_2$\big]}

\vspace{0.1cm}

\quad $Z_{j+1}:={\nu}[{\cal O}](Z_j)$ 
}}
\end{overview}

For some pseudo instructions, we will also give all instructions that are to be executed if such a pseudo instruction is to be performed. While {\sf instruction}$_1$, {\sf instruction}$_2$, $\ldots$ denote single instructions, each of the strings {\sf inst}, {\sf inst}$_1$, $\ldots$ may stand for a whole sequence of instructions that can be a single {\sf instruction} or a {\sf subprogram} or a {\sf pseudo instruction}. We sometimes also describe subprograms by sequences such as ${\sf inst}_{i_1}; {\sf inst}_{i_1+1}$ where we omit the labels, and so on. 

\vspace{0.3cm}

\noindent {\bf The pseudo instruction 
{\sf goto} $\ell$.}

Any command 

\fbox{{\sf goto} $\ell$} 

\vspace{0.1cm}

\noindent stands for the instruction 

\vspace{0.2cm}

\hspace{1.1cm}{\sf if $I_{1,1}=I_{1,1}$ then goto $\ell $ else goto $\ell$}.

or

\hspace{1.1cm}{\sf if $I_{1}=I_{1}$ then goto $\ell $ else goto $\ell$}.

\vspace{0.3cm}

\noindent {\bf The pseudo instructions $(Z_{2,1},Z_{2,2},\ldots):= (Z_{1,1},\ldots)$.}

The assignment

\vspace{0.2cm}

\fbox{$(Z_{2,1},\ldots, Z_{2,I_{2,1}}):= (Z_{1,1},\ldots, Z_{1,I_{1,1}})$}

\vspace{0.2cm}

\noindent stands for a subprogram as given in the first five lines in Overview \ref{CopyInp} (where the index registers and the labels used are examples). In this case, let $c(I_{1,7})=1$ and $c(I_{2,1})=1$ hold at the beginning. Let $c(I_{1,1})$ be the length $n$ of the tuple whose components are to be copied and let $c(I_{1,7})$ be the index of the register $Z_{1,I_{1,7}}$ whose content $c(Z_{1,I_{1,7}})$ is currently to be copied. In the same loop pass, $c(I_{2,1})$ is the index of the register $Z_{2,I_{2,1}}$ that should receive the value $c(Z_{1,I_{1,7}})$. 

\begin{overview}
[Program segment ${\cal P}_{\rm init}$ for copying the initial values]\label{CopyInp} 

\hfill

\nopagebreak 

\noindent \fbox{\parbox{11.8cm}{

\vspace{0.1cm}
{\sf 

\quad $\textcolor{blue}{1}:\, Z_{2,I_{2,1}}:=Z_{1,I_{1,7}};$

\quad $2 :\,$ {\sf if $I_{1,1}=I_{1,7}$ then goto $6$ else goto $3;$} 

\quad $3: \, I_{1,7}:=I_{1,7}+1;$

\quad $4: \,I_{2,1}:=I_{2,1}+1;$

\quad $5 :\,$ {\sf goto $1;$}

\quad $6 :\,$ goto $\textcolor{blue}{1'};$ 

}}}
\end{overview}

\vspace{0.1cm}

\noindent 
If this loop construct has been successfully processed (and label 6 is reached), then the content of $I_{2,1}$ is also $n$ which means that $c(I_{2,1})=n$ holds. The label \textcolor{blue}{$1'$ belongs to ${\cal P}_{\cal M}'$} (see Overview \ref{SimM}) where it is possible to use that $c(I_{1,1})+1$ is first the second index of the $Z$-register $Z_{1,n+1}$ that belongs to the first tape and could --- in case of a non-deterministic machine --- contain the first guess.

\vspace{0.3cm}

\noindent {\bf The pseudo instructions \,{\sf for $I_{1,3}:=1,2,\ldots, \big[I_{1,4}\big]$ do \{inst\}}.} 

Let {\sf inst} be a sequence of instructions, none of which may lead to a change of the content of $I_{1,3}$ \big[or $I_{1,4}$\big]. Then, the command 

\vspace{0.2cm}

\fbox{\sf for $I_{1,3}:=1,2,\ldots$ do inst} 

\vspace{0.2cm}

\noindent stands for 

\hspace{1.6cm} $I_{1,3}:=1;$ 

\hspace{1.1cm}$\overline \ell:\,\,$ {\sf inst}; $I_{1,3}:=I_{1,3} + 1$; {\sf goto $\overline \ell $ }.

\vspace{0.2cm}

\noindent This loop construct can be exited only if a command of the form {\sf goto $\ell$} is performed. $\ell$ can also be the label of the stop instruction. In the second case, we usually start with $c(I_{1,3})\leq c(I_{1,4})$. Then,

\vspace{0.2cm}

\fbox{\sf for $I_{1,3}:=1,2,\ldots, I_{1,4}$ do inst$; \,\,\, \ell : {\sf inst}_{\ell}$ } 

\vspace{0.2cm}

\noindent stands for 

\hspace{1.6cm} $I_{1,3}:=1$; {\sf inst}; 
\, {\sf if $I_{1,3}=I_{1,4}$ then goto $\ell$ else goto $\overline \ell; $}

\hspace{1.1cm}$\overline \ell:\,\,$ $I_{1,3}:=I_{1,3} + 1$; {\sf inst}; 
\, {\sf if $I_{1,3}=I_{1,4}$ then goto $\ell$ else goto $\overline \ell; $}

\hspace{1.1cm}$\ell:\,\,{\sf inst}_{\ell}$

\vspace{0.2cm}

\noindent {\bf The pseudo instructions $ I_{1,4}:=ca_1^{\big[+\big]}(I_{1,3});$ $I_{1,5}:=ca_2^{\big[+\big]}(I_{1,3})$.}

Let us use the pseudo instruction $I_{1,10}:=ca(I_{1,4},I_{1,5})$ that can be computed by using subprograms listed in \cite[Overview 2, p.\,588]{GASS20}. Then, any subprogram 

\fbox{$I_{1,4}:=ca_1(I_{1,3}); \,\,\,I_{1,5}:=ca_2(I_{1,3}); \,\,\, \ell+1: {\sf inst}_{\ell+1}$} 

\vspace{0.1cm}

\noindent stands for 

{\sf 
\hspace{1.55cm} for $I_{1,4}:=1,2,\ldots, I_{1,3}$ do 

\hspace{2.1cm} \{ for $I_{1,5}:=1,2,\ldots,I_{1,3}$ do 

\hspace{2.6cm} \{ $I_{1,10}:=ca(I_{1,4},I_{1,5});$ if $I_{1,3}=I_{1,10}$ then goto $\ell +1$ \}
 \}$;$

\hspace{1.2cm}$\ell:$ goto $\ell$; \hfill \textcolor{blue}{\it (ps.instr\,2)}

\hspace{0.6cm}$\ell+1: {\sf inst}_{\ell+1}$}

 \vspace{0.1cm}

If $ca_1(c(I_{1,3}))=0$ or $ca_2(c(I_{1,3}))=0$ holds, then the execution of the program cannot be continued and loops forever (which follows from executing the instruction labeled by $\ell$).
\fbox{$I_{1,4}:=ca_1^+(I_{1,3}); \,\,\,I_{1,5}:=ca_2^+(I_{1,3}); \,\,\, \ell+1: {\sf inst}_{\ell+1}$} results from replacing \textcolor{blue}{\it (ps.instr\,2)} by

\vspace{0.1cm}

\hspace{1.2cm}$\ell:$ $I_{1,4}:=1$; $I_{1,5}:=1$;

\vspace{0.3cm}

\noindent {\bf The pseudo instruction {\sf if $I_{1,2}=i$ then goto $i\,'$}.} 

If the value $i$ of $I_{1,2}$ is the result of the execution of ${\sf inst}_k'$ including $I_{1,2}:=k+1$ considered in Overview \ref{SimM}, then $i$ is the label $k+1$ of the instruction that follows a $\nu$-instruction in ${\cal P}_{\cal M}$ and therefore $i>1$ holds. If, moreover, the following pseudo instruction labeled by $\widetilde 4$ is used to complete the process of evaluating a $\nu$-instruction as suggested in Overview \ref{SimNuInst} and $2', 3', \ldots, \ell_{{\cal P}_{\cal M}}'$ are the only labels in ${\cal P}_{\cal M}'$ derived directly from the labels $2, 3, \ldots, \ell_{{\cal P}_{\cal M}}$ in Overview \ref{SimM}, then

\vspace{0.2cm}

\fbox{$\widetilde 4:\,\,$ {\sf if $I_{1,2}=i$ then goto $i\,'$}}

\vspace{0.2cm}

\noindent stands for the following sequence of labeled subprograms.

\vspace{0.1cm}

\hspace{1cm }$\widetilde 4:\,\,$ {\sf if $I_{1,2}=2$ then goto $2\,'$ else goto $\widetilde 5$};

\hspace{1cm }$\ldots$ 

$\widetilde {\ell_{{\cal P}_{\cal M}}\!+1}:\,\,$ {\sf if $I_{1,2}=\ell_{{\cal P}_{\cal M}}\!-1$ then goto $(\ell_{{\cal P}_{\cal M}}\!-1)'$ else goto $\widetilde {\ell_{{\cal P}_{\cal M}}+2}$};

$\widetilde {\ell_{{\cal P}_{\cal M}}\!+2}:\,\,$ {\sf if $I_{1,2}=\ell_{{\cal P}_{\cal M}}$ then goto $\ell_{{\cal P}_{\cal M}}'$ }

\vspace{0.2cm}

\noindent These subprograms can be replaced by permitted instructions as follows.

\begin{overview}[Determining the next label depending on an index]\label{EvaluForLabel}

\hfill

\nopagebreak 

\noindent \fbox{\parbox{11.8cm}{

\vspace{0.1cm}

\hspace{1cm }$\widetilde 4:\,\,$ $I_{1,8}:=1;$ 

\hspace{1.5cm } $I_{1,8}:=I_{1,8}+ 1;$ {\sf if $I_{1,2}=I_{1,8}$ then goto $2\,'$ else goto $\widetilde 5;$}

\hspace{1cm }$\widetilde 5:\,\,$ $I_{1,8}:=I_{1,8}+ 1;$ {\sf if $I_{1,2}=I_{1,8}$ then goto $3\,'$ else goto $\widetilde 6;$}

\hspace{1cm }$\ldots$ 

\,\, $\widetilde {\ell_{{\cal P}_{\cal M}}\!+1}:\,\,$ $I_{1,8}:=I_{1,8}+ 1;$ {\sf if $I_{1,2}=I_{1,8}$ then goto $(\ell_{{\cal P}_{\cal M}}\!-1)'$ else goto $\widetilde {\ell_{{\cal P}_{\cal M}}\!+2};$}

\,\, $\widetilde {\ell_{{\cal P}_{\cal M}}\!+2}:\,\,$ $I_{1,8}:=I_{1,8}+ 1;$ {\sf if $I_{1,2}=I_{1,8}$ then goto $\ell_{{\cal P}_{\cal M}}'$ else goto $\ell_{{\cal P}_{\cal M}}'$}}}
\end{overview}

\begin{figure}
{\scriptsize
\hspace{6.6cm}$I_{k_{\cal M}+1}$

\hspace{7cm}$\downarrow$

(1)\hspace{2cm}\begin{tabular}{|c|c|c|c|c|c|c|c|c|c|c|c|}\hline
$z_1$&$z_2$&$\ldots$&$z_i$&$\ldots$&$z_{k}$&$x_n\!$&$u_1$&$u_2$&$u_3$&$u_4$&$\ldots$ \\\hline
\end{tabular}

\hspace{5.3cm}$\uparrow$ \hspace{1.47cm}$\uparrow$

\hspace{5.3cm}$I_j$ \hspace{0.9cm} $I_{k_{\cal M}+2}$

$\,\,\downarrow$ $\ldots; \textcolor{blue}{I_{k_{\cal M}+1}:=I_{k_{\cal M}+1}+1}$ 

\hspace{7.33cm}$I_{k_{\cal M}+1}$

\hspace{7.7cm}$\downarrow$

(2)\hspace{2cm}\begin{tabular}{|c|c|c|c|c|c|c|c|c|c|c|c|}\hline
$z_1$&$z_2$&$\ldots$&$z_i$&$\ldots$&$z_{k}$&$x_n\!$&$u_1$&$u_2$&$u_3$&$u_4$&$\ldots$ \\\hline
\end{tabular}

\hspace{5.6cm}$\uparrow$ \hspace{1.17cm}$\uparrow$

\hspace{5.6cm}$I_j$ \hspace{0.6cm} $I_{k_{\cal M}+2}$

$\,\,\downarrow\,$ \textcolor{blue}{$Z_{I_{k_{\cal M}+1}}:=Z_{I_{k_{\cal M}+2}}$}

\hspace{7.33cm}$I_{k_{\cal M}+1}$

\hspace{7.7cm}$\downarrow$

(3)\hspace{2cm}\begin{tabular}{|c|c|c|c|c|c|c|c|c|c|c|c|}\hline
$z_1$&$z_2$&$\ldots$&$z_i$&$\ldots$&$z_{k}$&$x_n\!$&$x_n\!$&$u_2$&$u_3$&$u_4$&$\ldots$ \\\hline
\end{tabular}

\hspace{5.6cm}$\uparrow$ \hspace{1.17cm}$\uparrow$

\hspace{5.6cm}$I_j$ \hspace{0.6cm} $I_{k_{\cal M}+2}$

 $\,\,\downarrow\,$ \textcolor{blue}{$I_{k_{\cal M}+2}:=I_{k_{\cal M}+2}+1$}

\hspace{7.33cm}$I_{k_{\cal M}+1}$

\hspace{7.7cm}$\downarrow$

(4)\hspace{2cm}\begin{tabular}{|c|c|c|c|c|c|c|c|c|c|c|c|}\hline
$z_1$&$z_2$&$\ldots$&$z_i$&$\ldots$&$z_{k}$&$x_n\!$&$x_n\!$&$u_2$&$u_3$&$u_4$&$\ldots$ \\\hline
\end{tabular}

\hspace{5.6cm}$\uparrow$ \hspace{1.83cm}$\uparrow$

\hspace{5.6cm}$I_j$ \hspace{1.3cm} $I_{k_{\cal M}+2}$}

\caption{The execution of ${\sf init}(Z_{I_{k_{\cal M}+1}})$}\label{InitZReg}

\end{figure}

\vspace{0.3cm}

\noindent {\bf The pseudo instructions ${\sf init}(Z_{d_1,I_{d_1,k_{\cal M}+1}})$ and ${\sf init}(Z_{I_{k_{\cal M}+1}})$.}

The pseudo instruction 

\fbox{${\sf init}(Z_{d_1, I_{d_1,k_{\cal M}+1}})$} 

\vspace{0.1cm}

\noindent stands for 

\hspace{1.1cm} $ I_{d_1,k_{\cal M}+1}:=I_{d_1,k_{\cal M}+1}+1$; 
\,$Z_{d_1,I_{d_1,k_{\cal M}+1}}:=Z_{d_1,I_{d_1,k_{\cal M}+2}}$;

\hspace{1.1cm} $I_{d_1,k_{\cal M}+2}:=I_{d_1,k_{\cal M}+2}+1$

\vspace{0.2cm}

\noindent The execution of such an additional subprogram in simulating the execution of an instruction of type (7) can be important to be prepared for the case that the $Z$-registers $Z_{d_1,1}, Z_{d_1,2},\ldots$ of a machine ${\cal N}$ must be used (instead of $Z_1,Z_2,\ldots$) in simulating a BSS RAM ${\cal M}$ in ${\sf M}_{\cal A}$ or ${\sf M}_{\cal A}^{\nu}(Q)$ for some oracle $Q$. Then, it is possible that a $Z$-register $Z_s$ ($s\geq n$) of ${\cal M}$ and the corresponding $Z$-register $Z_{d_1,s}$ of ${\cal N}$ are to be considered (again). If $\vec x$ is the input for ${\cal M}$, then, at the beginning of a simulation, it must be ensured that the register $Z_s$ of ${\cal M}$ contains the value $x_n$ and the same should therefore hold for $Z_{d_1,s}$. Consequently, to be on the safe side, $x_n$ is copied into $Z_{d_1,s}$ before it could possibly be used. 

Figure \ref{InitZReg} illustrates this process for \fbox{${\sf init}(Z_{I_{k_{\cal M}+1}})$}. In (1) of Figure \ref{InitZReg}, the $Z$-registers contain the individuals $z_1, z_2, \ldots, z_{i}$, $ \ldots, z_{k}, x_n, u_1, u_2, u_3, u_4,\ldots$. An arrow at the top or bottom of the box of a $Z$-register indicates which of the index registers $I_{k_{\cal M}+1}$ and $I_{k_{\cal M}+2}$ contains the index or the address of this $Z$-register.

\vspace{0.3cm}

\noindent {\bf The pseudo instruction ${\sf initguess}(Z_{I_{k_{\cal M}+1}})$.}

This pseudo instruction 

\fbox{${\sf initguess}(Z_{I_{k_{\cal M}+1}})$} 

\vspace{0.1cm}

\noindent stands for 

\hspace{1.1cm} $I_{k_{\cal M}+1}:=I_{k_{\cal M}+1}+1$; \, $I_{k_{\cal M}+2}:=I_{k_{\cal M}+1}+1$; \, $I_{k_{\cal M}+3}:=1$; 

\hspace{1.1cm} $Z_{I_{k_{\cal M}+2}}:=Z_{I_{k_{\cal M}+3}}$; \, $Z_1:=\nu({\cal O})(Z_{1},\ldots,Z_{I_{1}})$;\, 

\hspace{1.1cm} $Z_{I_{k_{\cal M}+1}}:=Z_{I_{k_{\cal M}+3}}$; \, $Z_{I_{k_{\cal M}+3}}:=Z_{I_{k_{\cal M}+2}}$

\vspace{0.2cm}

\noindent This subprogram enables the simulation of a non-deterministic machine by using $U_{\cal A}^\infty$ as oracle without carrying out the guessing process as part of the input procedure as in the case of a non-deterministic BSS RAM in ${\sf M}_{\cal A}^{\rm ND}$.

\subsection{The non-deterministic simulation of $\nu$-oracle machines}

Let ${\cal A}\in {\sf Struc}$, $Q$ be a set in ${\rm SDEC}_{\cal A}$, and ${\cal N}_Q$ be one of the BSS RAMs in ${\sf M}_{\cal A}$ that semi-decide $Q$. We will give more details why $({\rm SDEC}_{\cal A}^\nu)^Q\subseteq {\rm SDEC}_{\cal A}^{\rm ND}$ holds (cf.\,\,Theorem \ref{OracleMNondetM}). Therefore, let ${\cal M}$ be any BSS RAM in ${\sf M}_{\cal A}^{\nu}(Q)$. We want to define a suitable non-deterministic machine ${\cal N}_{\cal M}\in ({\sf M}^{(3)}_{\cal A})^{\rm ND}$ such that ${\cal N}_{\cal M}$ can simulate ${\cal M}$ on any input $\vec x\in U_{\cal A}^\infty$ which implies $H_{{\cal N}_{\cal M}}= H_{\cal M}$. Roughly speaking, each $\nu$-instruction in ${\cal P}_{\cal M}$ must be replaced by a subprogram containing a suitably adapted copy of ${\cal P}_{{\cal N}_Q}$ such that ${\cal N}_{\cal M}$ can also simulate ${\cal N}_Q$ for inputs resulting from the use of its own guesses (produced by ${\rm Input}_{{\cal N}_{\cal M}}$) as components in addition to the input values of the considered $\nu$-instruction. For this purpose, we define some program segments where two of which, ${\cal P}_{\cal M}'$ and ${\cal P}_{{\cal N}_Q}''$, are derived from the programs ${\cal P}_{\cal M}$ and ${\cal P}_{{\cal N}_Q}$ (cf.\,\,Overviews \ref{SimM} and \ref{SubprSimPNQ}). 
For the simulation of ${\cal M}$, let ${\cal N}_{\cal M}$ be equipped with the following 3 tapes. 

\vspace{0.1cm}

 {\bf Tape 1.} The first tape $Z_{1,1},Z_{1,2},\ldots$ should contain the input $\vec x$ and the guesses of ${\cal N}_{\cal M}$.
We assume that any $n$-ary input $\vec x$ for ${\cal N}_{\cal M}$ stored in $Z_{1,1}$, \ldots, $Z_{1,n}$ is also the input for ${\cal M}$ whose work is to be evaluated. When starting the simulation, $m_0$ should be $0$. Consequently, $Z_{1,n+1}, Z_{1,n+2},\ldots$ contain the guesses $y_{m_0+1},y_{m_0+2},\ldots$ of ${\cal N}_{\cal M}$ (for $m_0=0$) some (or all) of which can be used if the first execution of a $\nu$-instruction has to be simulated. 

\vspace{0.1cm}

{\bf Tape 2.} The second tape $Z_{2,1},Z_{2,2},\ldots$ should contain the values $c(Z_1)$, $c(Z_2)$, $\ldots$ of the $Z$-registers $Z_1,Z_2,\ldots$ of ${\cal M}$ for executing the instructions, that correspond to the instructions of types (1) to (8) in ${\cal P}_{\cal M}$, in the same way as ${\cal M}$. Before ${\cal N}_{\cal M}$ simulates ${\cal M}$, the input values $x_1,\ldots,x_n$ of ${\cal N}_{\cal M}$ stored in $Z_{1,1}$, \ldots, $Z_{1,n}$ are copied from the first tape into the first $Z$-registers $Z_{2,1}$, \ldots, $Z_{2,n}$ of the second tape and $I_{2,1}$ receives the value of $ I_{1,1}$. 

\vspace{0.1cm}

{\bf Tape 3.} The third tape $Z_{3,1},Z_{3,2},\ldots$ should be used to evaluate the $\nu$-instructions of ${\cal M}$ and simulate ${\cal N}_Q$ in this process. If $c(I_{2,1})=n_0$ and the execution of a $\nu$-instruction should be simulated, then the first $n_0$ values $z_1,\ldots, z_{n_0}$ of the $Z$-registers of ${\cal M}$ are copied from $Z_{2,1},\ldots, Z_{2,n_0}$ into $Z_{3,1},\ldots, Z_{3,n_0}$ and used in simulating ${\cal N}_Q$ on $(c(Z_{2,1}),\ldots, c(Z_{2,n_0}), c(Z_{1,n+m_0+1}), \ldots, c(Z_{1,n+m_0+m}))$ for certain $m_0\geq 0$ and suitable integers $m\geq 1$. 

\begin{overview}
[Program ${\cal P}_{{\cal N}_{\cal M}}$ of ${\cal N}_{\cal M}$]

\hfill

\nopagebreak 

\noindent \fbox{\parbox{11.8cm}{

\vspace{0.1cm} 

\quad ${\cal P}_{\rm init}$ \hfill{\scriptsize(The {\sf first label is \textcolor{blue}{$1$}}.)}

\quad ${\cal P}_\nu$  \hfill{\scriptsize(Its first label is \textcolor{blue}{$1^*$}, its last label is $\widetilde 4$.)}

\quad ${\cal P}_{\cal M}'$ \hfill{\scriptsize(Its first label is \textcolor{blue}{$1'$},  the {\sf  last label is $\textcolor{blue}{\ell _{{\cal P}_{\cal M}}'}=\ell_{{\cal P}_{{\cal N}_{\cal M}}}$}.)}

}}
\end{overview}
Before ${\cal N}_{\cal M}$ starts its program, its first tape receives the input $\vec x$ and guessed values and all registers of its second tape receive the value $x_n$ by applying its multi-valued input procedure. For simulating ${\cal M}$ by ${\cal N}_{\cal M}$, the own input is copied onto the second tape of ${\cal N}_{\cal M}$ by means of the program segment ${\cal P}_{\rm init}$ given in Overview \ref{CopyInp}. As a consequence of executing ${\cal P}_{\rm init}$, the second tape contains --- before ${\cal N}_{\cal M}$ begins the simulation of ${\cal M}$ on its input $\vec x$ by means of ${\cal P}_{\cal M}'$ --- the values $x_1,\ldots, x_n,x_n, x_n,\ldots$ which correspond to the start configuration of ${\cal M}$ on input $\vec x$. Overview \ref{SimM} gives the details how ${\cal P}_{\cal M}'$ can be derived from ${\cal P}_{\cal M}$. 

\begin{overview}
[Program segment ${\cal P}_{\cal M}'$ for simulating $\cal M\in{\sf M}_{\cal A}^{\nu}(Q)$]\label{SimM}

\hfill

\nopagebreak 

\noindent \fbox{\parbox{11.8cm}{

Let the program ${\cal P}_{\cal M}$ of the 1-tape machine ${\cal M}\in{\sf M}_{\cal A}^{\nu}(Q)$ be of the form \textcolor{blue}{$(*)$} (see p.\,\,\pageref{star}) and let ${\cal P}_{\cal M}'$ belong to ${\cal P}_{{\cal N}_{\cal M}}$ and given by

\vspace{0.2cm}

$\quad\textcolor{blue}{1\,'}\, : {\sf inst} '_1 ; \quad\ldots;\,\, ( \ell _{{\cal P}_{\cal M}}-1)\,' : {\sf inst} '_{ \ell _{{\cal P}_{\cal M}}-1} ;\,\, \textcolor{blue}{\ell _{{\cal P}_{\cal M}} '} : {\sf stop}.$

\vspace{0.2cm}

where each {\sf inst}$'_k$ results from {\sf instruction}$_k$ in ${\cal P}_{\cal M}$ by replacing
\begin{itemize} 
\item\parskip -0.5mm 
 each symbol $Z_i$ by $Z_{2,i}$, 
\item each symbol $I_i$ by $I_{2,i}$, and 
\item each command \,{\sf goto $\ell$}\, by \,{\sf goto $\ell\,'$}
\end{itemize}\setlength{\parskip}{-0.5em}
 if \,{\sf instruction}$_k$\, is not of type (10), or otherwise it is determined as follows.}}
\nopagebreak 

\noindent \fbox{\parbox{11.8cm}{

If \,{\sf instruction}$_k$\, has the form \fbox{$Z_j:=\nu[{\cal O}](Z_1,\ldots,Z_{I_1})$}, then let \,${\sf inst}'_k$\, be the subprogram 

\qquad \,\,\,\, $I_{1,2}:= k +1; $ \, $I_{2,k_{\cal M}+1}:=j;$ \, {\sf goto} \textcolor{blue}{$1^*$} \,.
}}\end{overview}
The subprogram labeled by \textcolor{blue}{$1^*$ belongs to ${\cal P}_\nu$} and is given in Overview \ref{SimNuInst}. It allows to copy the values $c(Z_{2,1}),\ldots,c(Z_{2,c(I_{2,1})})$ onto tape 3. After this, an implementation of the algorithm {\it (algo\,1)} whose application was proposed in the proof of Theorem \ref{OracleMNondetM} follows. For details, see Section \ref{SimOracleInstr} where the following index registers are also used.

\begin{overview}
[The use of index registers $I_{1,1},I_{1,2},\ldots$ of ${\cal N}_{\cal M}$]\label{SimNuInstValues}

\hfill

\nopagebreak 

\noindent \fbox{\parbox{11.8cm}{

\setlength{\extrarowheight}{-2pt} 

\begin{tabular}{lll}
& {\rm Meaning}&{\rm Values} \\
$I_{1,1}\!$&index of the last used individual on tape 1\\
 &after evaluating the $i_0^{\rm th}$ $\nu$-instruction &$n +\sum_{i=1}^{i_0}m_i $\\
$I_{1,2}\!$&contains a value for determining a label $i'$&$i$ \\
$I_{1,3}\!$&auxiliary index for counting the loops&$s$ \\
$I_{1,4}\!$& determined by executing $I_{1,4}:=ca_1^+(I_{1,3}) $ &$m$ \\
$I_{1,5}\!$& determined by executing $I_{1,5}:=ca_2^+(I_{1,3})\big[\!\!+\!1\big]$&$s_0$ and $s_0+1$\\
$I_{1,6}\!$&index of the last input value on tape 3\\
 &for evaluating one execution of $s_0$ steps of ${\cal N}_Q\!\!$\\
 &determined by executing $I_{1,6}:= I_{3,1}+I_{1,4}$&$n_0+m$ \\
$I_{1,7}\!$&auxiliary index used in ${\cal P}_{\rm init}$ and ${\cal P}_{\rm init}^*$& $1\,..\,n \!+\!\sum_{i=1}^{i_0}m_i \!+\!m$\\
$I_{1,8}\!$&auxiliary index for evaluating $I_{1,2}$&$1\,.. \, 2\,..\,\ell_{{\cal P}_{\cal M}}$\\
$I_{1,9}\!$&auxiliary index for counting the steps &$1\,..\,s_0+1$\\
$\!\!I_{1,10}\!$& auxiliary index for computing $m$ and $s_0$ from $s\!\!\!$&$1\,..\, ca(s,s)\!\!$\\\end{tabular}
}}
\end{overview}

\subsection{Replacing $\nu$-instructions by trying to semi-decide $Q$}\label{SimOracleInstr}
The subprogram of ${\cal P}_\nu$ whose first instruction is labeled by $\widetilde 1$ (given in Overview \ref{SimNuInst}) shows how \fbox{$k:\,\, Z_j:=\nu[{\cal O}](Z_1,\ldots,Z_{I_1})$} can be replaced so that the execution of such an oracle instruction can be simulated by the repeated use of a suitably adapted copy ${\cal P}_{{\cal N}_Q}''$ of ${\cal P}_{{\cal N}_Q}$ (given in Overview \ref{SubprSimPNQ}). 

{\bf The used guesses.} Due to its input procedure, the guesses of ${\cal N}_{\cal M}$ are stored in $Z_{1,n+1}, Z_{1,n+2},\ldots$. In {\it (algo\,1)}, we have used $m_0$ for counting the guesses processed in the simulation of the first $i_0$ executions of $\nu$-instructions. According to this, for any $\nu$-instruction in the sequence of instructions performed by ${\cal M}$ (whose execution is to be simulated), there is a new $m_0$ resulting from $m_0:=0$ (for $i_0=0$ and the evaluation of the first $\nu$-instruction), $m_0:=m_1$ (for $i_0=1$), $m_0:=m_1+m_2$ (for $i_0=2$),\ldots. Depending on this, ${\cal N}_{\cal M}$ can use some of its guesses: (G0) $y_1,\ldots, y_{m_1}$ (for the first $\nu$-instruction), (G1) $y_{m_1+1},\ldots, y_{m_1+m_2}$ (for the second $\nu$-instruction), (G2) $y_{m_1+m_2+1},\ldots$, $y_{m_1+m_2+m_3}$, $\ldots $. If, for $i_0\in \{0,1,\ldots\}$, the simulations of the first $i_0$ executions of certain $\nu$-instructions were successfully completed, then the label $\ell_{{\cal P}_{{\cal N}_Q}}''$ corresponding to $\ell_{{\cal P}_{{\cal N}_Q}}$ was reached $i_0$ times and, for $i_0>0$, each of the individuals in $\{y_1,\ldots, y_{m_0}\}$ (with $m_0={\sum_{i=1}^{i_0} m_i}$) were already used and $y_1, y_{m_1+1}, \ldots, y_{\sum_{i=1}^{i_0-1} m_i+1}$ are (or were) the results of $\nu$-instructions copied onto the second tape. After this, 

\vspace{0.2cm}

$({\rm G}i_0)\qquad y_{\sum_{i=1}^{i_0} m_i+1}, y_{\sum_{i=1}^{i_0} m_i+2}, \ldots$\hfill (intended for the $(i_0+1)^{\rm th}$ execution)

\vspace{0.2cm}

\noindent can enter into the evaluation of the $(i_0+1)^{\rm th}$ $\nu$-instruction by ${\cal N}_{\cal M}$ if it is necessary to simulate the execution of a further $\nu$-instruction. Any $\nu$-instruction itself can be evaluated by performing the algorithm given in Overview \ref{SimNuInst}.

We would like to point out that $i_0$ and $m_0$ are used only in the description of the simulation and that they are not directly stored in index registers that are changed by applying the program. For describing the simulation of the execution of any further $\nu$-instruction (by executing the subprogram labeled by $\widetilde 1$ to $\widetilde 2$ in ${\cal P}_\nu$), the integer $i_0$ must be incremented by 1 (which means that we apply $i_0:=i_0+1$ or $i_0$++) and the integer $m_0$ (also used only in a less formal description) obtains a new value, i.e., the value $\sum_{i=1}^{i_0}m_i$ by adding $m_{i_0}$, given by $m_{i_0}:=m$, to the old value $m_0$ such that, after executing the instruction labeled by $\widetilde 3$, we have $m_0=c(I_{1,1})-n$.

{\bf The repeated use of ${\cal P}_{{\cal N}_Q}''$.} The execution of each $\nu$-instruction itself can be simulated by using a loop construct as summarized in {\it (algo\,1)} since, at the beginning, it is not clear how many guesses and how many steps are necessary for reaching $\ell_{{\cal N}_Q}''$. For this purpose, ${\cal N}_{\cal M}$ repeats the execution of the loop block for $s= 1,2,\ldots$ (and $m=ca_1^+(s)$) until the label $\ell_{{\cal P}_{{\cal N}_Q}}''$ is reached. This means that, for each considered $s$, ${\cal N}_{\cal M}$ starts the simulation of ${\cal N}_Q$ again and simulates $ca_2^+(s)$ steps of ${\cal N}_Q$ on the initial values whose number is determined by $ca_1^+(s)$ and, consequently, given by $c(I_{3,1})=c(I_{1,6})=n_0+m$ (computed by using subprograms labeled by $\overline 1$ and $3^*$ in Overviews \ref{SimNuInst} and \ref{CopyInp2}, respectively) and

\vspace{0.3cm}

$({\rm Inp})\hspace{2cm}\begin{array}{c}
c(Z_{3,1})= z_1, \ldots, c(Z_{3,n_0})= z_{n_0},\\c(Z_{3,n_0+1})= y_{m_0+1}, \ldots, c(Z_{3,n_0+m})= y_{m_0+m}.\end{array}$

\vspace{0.3cm}

\noindent Let us assume that the registers of ${\cal N}_{\cal M}$ contain the following values before ${\cal N}_{\cal M}$ simulates the $(i_0+1)^{\rm th}$ execution of a $\nu$-instruction of ${\cal P}_{\cal M}$. 

\begin{overview}[Values of registers before evaluating a $\nu$-instruction]

\hfill

\nopagebreak 

\noindent \fbox{\parbox{11.8cm}{

\begin{itemize}
\item\parskip -0.7mm 
$n_0=c(I_{2,1}) $ is also the content of the registers
 $I_1$ of ${\cal M}$,
\item $(z_1,\ldots, z_{n_0})$ is stored in 
\begin{itemize}
\item\parskip -0.7mm 
$Z_{2,1},\ldots, Z_{2,n_0}$ of ${\cal N}_{\cal M}$,
\end{itemize}
and the content of $Z_1,\ldots, Z_{n_0}$ of ${\cal M}$,
\item for $m_0=\sum_{i=1}^{i_0}m_i$, the value $n+m_0$ is stored in $I_{1,1}$ of ${\cal N}_{\cal M}$.
\end{itemize}
}}
\end{overview}
If the input of a $\nu$-instruction of ${\cal M}$ is $(z_1,\ldots, z_{n_0})$ and stored in $Z_{3,1},\ldots,$ $Z_{3,c(I_{3,1})}$ (after executing the subprogram labeled by $1^*$ in Overview \ref{SimNuInst}), then ${\cal N}_{\cal M}$ simulates the $(i_0+1)^{\rm th}$ execution of a $\nu$-instruction by using the subprogram given in the lines 2 to 6 of Overview \ref{SimNuInst}.

\begin{overview}
[Program segment ${\cal P}_\nu$ for processing a $\nu$-instruction]\label{SimNuInst}

\hfill

\nopagebreak 

\noindent \fbox{\parbox{11.8cm}{

\vspace{0.1cm} 

\hspace{0.08cm} $\textcolor{blue}{1^*}: \, (Z_{3,1},\ldots, Z_{3,I_{3,1}}):= (Z_{2,1},\ldots, Z_{2,I_{2,1}});$ \hfill {\scriptsize (Recall, $c(I_{2,1})=n_0$.)}

\quad{\sf $\widetilde 1:\,$ $I_{1,3}:=1;$ 

\quad $\overline 1:\,$ $I_{1,4}:=ca_1^+(I_{1,3}) ;\,\, I_{1,5}:=ca_2^+(I_{1,3}) ;\,\, I_{1,6}:= I_{3,1}+I_{1,4};\,I_{1,9}:=1;$

\hspace{0.1cm} ${\cal P}^*_{\rm init}$

\hspace{0.1cm} ${\cal P}_{{\cal N}_Q}''$} \hfill {\scriptsize (i.e., simulate $c(I_{1,5})$ steps of ${\cal N}_Q$ on $(c(Z_{3,1}),\ldots, \underbrace{c(Z_{3,c(I_{3,1})+1}),\ldots,c(Z_{3,c(I_{1,6})}}_{{\it determined \,\,by\,\,}{\cal P}^*_{\rm init} }))$)}

{\sf
\quad $ \textcolor{blue}{\overline{2}}: \,\,$ $I_{1,3}:=I_{1,3}+1;$\,\, goto $\bar 1;$}

\quad $\textcolor{blue}{\widetilde 2}:\,\,$ $Z_{2,I_{2,k_{\cal M}+1}}:=Z_{1,I_{1,1}+1};$
\hfill {\scriptsize(Recall, $c(Z_{1,c(I_{1,1})+1})= y_{m_0+1}$.)}

\quad $\widetilde 3:\,\,$ $I_{1,1}:=I_{1,1} + I_{1,4};$ 
\hfill {\scriptsize (Recall, $c(I_{1,4})= m$.)}

\quad $\widetilde 4:\,\,$ {\sf if $I_{1,2}=i$ then goto $\textcolor{blue}{i\,'};$} 
\hfill {\scriptsize(Recall,  \textcolor{blue}{$i\,'$ belongs to ${\cal P}_{\cal M}'$}, cf.\,\,Overview \ref{EvaluForLabel}.)}
}}\end{overview}

When starting a simulation of $ca_2^+(s)$ steps of ${\cal N}_Q$ during the evaluation of the $(i_0+1)^{\rm th}$ $\nu$-instruction, the following values are the initial values.

\begin{overview}[Values of registers before simulating $s_0$ steps of ${\cal N}_Q$]
\hfill

\nopagebreak 

\noindent \fbox{\parbox{11.8cm}{

After executing the subprograms labeled by $1^*$, $\widetilde 1$, and $\overline 1$ in Overview \ref{SimNuInst} and ${\cal P}_{\rm init}^*$ given in Overview \ref{CopyInp2},
 the following values are stored in registers of ${\cal N}_{\cal M}$:

\begin{itemize}
\item\parskip -0.5mm 
$s=c(I_{1,3})$,
\quad $s_0=c(I_{1,5})$,
\quad $m=c(I_{1,4})$,
\quad $n_0=c(I_{3,1})$,
\item 
 $(z_1,\ldots, z_{n_0})$ is stored in $Z_{3,1},\ldots, Z_{3,n_0}$,
\item
 $(y_{\sum_{i=1}^{i_0}m_i+1},\ldots, y_{\sum_{i=1}^{i_0}m_i+m})$ is stored in $Z_{3,n_0+1},\ldots, Z_{3,n_0+m}$.
\end{itemize}
}}
\end{overview}
If ${\cal N}_Q$ halts on $(z_1,\ldots, z_{n_0},y_{m_0+1},\ldots, y_{m_0+m}) $, then $y_{m_0+1}$ can be transferred by using the instruction labeled by $\widetilde 2$. The result is $c(Z_{2,j})=y_{m_0+1}$. Since ${\cal N}_{\cal M}$ needs $m_0:= m_0+m$ for the execution of a further $\nu$-instruction (if necessary), ${\cal N}_{\cal M}$ must execute the subprogram labeled by $\widetilde 3$. The label $i'$ of the next instruction in ${\cal P}_{\cal M}'$ that is to be executed is determined by the pseudo instruction labeled by $\widetilde 4$.

\begin{overview}[Values of registers after evaluating a $\nu$-instruction]
\hfill

\nopagebreak 

\noindent \fbox{\parbox{11.8cm}{

When label $\,\widetilde 2$ is reached:

$y_{m_0+1} =c(Z_{3,n_0+1})= c(Z_{1,c(I_{1,1})+1})$ 

\begin{itemize}
\item\parskip -0.7mm 
is an output value of the $\nu$-instruction whose execution was simulated,
\item will be assigned to $Z_{2,j}$ (corresponding to $Z_j$ of ${\cal M}$).
\end{itemize}
}}
\end{overview}

{\bf The repeated use of the same guesses.}
During the evaluation of the $(i_0+1)^{\rm th}$ $\nu$-instruction by executing ${\cal P}_\nu$, ${\cal N}_{\cal M}$ uses the same sequence of guesses considered in (G$i_0$) and stored in $Z_{1,n+m_0+1}$, $Z_{1,n+m_0+2}, \ldots$ for all simulations of ${\cal N}_Q$ on its inputs $(z_1,\ldots, z_{n_0},y_{m_0+1},\ldots)$ --- that are stored in $c(Z_{3,1}),\ldots,c(Z_{3,n_0+m})$ as described by (Inp) --- for integers $m$ increasing from time to time. If there are integers $s$ for which ${\cal N}_Q$ halts on an input given by (Inp) for $m=ca_1^+(s)$ after $ca_2^+(s)$ steps, then ${\cal N}_{\cal M}$ can reach the label $\ell_{{\cal P}_{{\cal N}_Q}}''$. For the smallest $s$ of this kind, ${\cal N}_{\cal M}$ finishes the execution of the loop construct and thus the simulation of the execution of the corresponding $\nu$-instruction by executing {\sf goto $\widetilde 2$}. Otherwise, if $\ell_{{\cal P}_{{\cal N}_Q}}''$ cannot be reached, then the execution of the loop construct does not terminate and ${\cal N}_{\cal M}$ uses more and more guesses that are stored in $Z_{1, n+\sum_{i=1}^{i_0} m_i+1}, Z_{1, n+\sum_{i=1}^{i_0} m_i+2},\ldots$ and that are copied, for each $s$, again into $Z_{3, n_0+1}, Z_{3, n_0+2},\ldots$ at the beginning of performing the instructions in the loop construct. 

{\bf The use of further guesses.} If the simulation comes to the instruction labeled by $\widetilde 3$ after executing the instruction labeled by $\widetilde 2$, then we can additionally use $i_0:=i_0+1$, $m_{i_0}=m$, and $m_0=\sum_{i=1}^{i_0}m_i$ in the informal description of the program ${\cal N}_{\cal M}$. That helps to express that new guesses $y_{\sum_{i=1}^{i_0} m_i+1}$, $y_{\sum_{i=1}^{i_0} m_i+2}, \ldots, $ are used as new additional input components for simulating the execution of the next (which means the new $(i_0+1)^{\rm th}$) $\nu$-instruction. If ${\cal N}_{\cal M}$ would need more than its own guesses in the attempt to simulate ${\cal M}$, then ${\cal N}_{\cal M}$ uses also the last component of its own input. If ${\cal M}$ executes its stop instruction after $t$ steps, then its simulation requires at most $\sum_{i=0}^t m_i$ guesses of ${\cal N}_{\cal M}$ (which can be provided by the multi-valued input procedure of ${\cal N}_{\cal M}$).

\subsection{More details: The simulation of ${\cal N}_Q$ semi-deciding $Q$}

 For evaluating a $\nu$-instruction of ${\cal P}_{\cal M}$ for the input $(z_1,\ldots, z_{n_0})$ on tape 3 by ${\cal N}_{\cal M}$, let $I_{1,3}$ be the register for storing the loop counter $s$ and let $c(I_{1,4})=ca_1^+(c(I_{1,3}))$ and $c(I_{1,5})=ca_2^+(c(I_{1,3}))$ be determined by pseudo instructions (see Overview \ref{PseudoInst}) whose execution implies $m=c(I_{1,4})$ and $s_0=c(I_{1,5})$. The program segment ${\cal P}_{\rm init}^*$ given in Overview \ref{CopyInp2} allows to simulate the application of the input procedure of ${\cal N}_Q\in{\sf M}_{\cal A}$ as far as possible. By executing the pseudo instruction labeled by $1^*$ in Overview \ref{SimNuInst}, the registers $Z_{3,1},\ldots, Z_{3,n_0}$ obtain the values $z_1,\ldots,z_{n_0}$ that are stored in $Z_{2,1},\ldots, Z_{2,n_0}$. By executing the pseudo instruction labeled by $2^*$, $y_{m_0+1}$, \ldots, $y_{m_0+m} $ stored in $Z_{1, n+ \sum_{i=1}^{i_0} m_i+1}, \ldots, Z_{1,n+\sum_{i=1}^{i_0} m_i +m}$ are copied into $Z_{3,n_0+1},\ldots, Z_{3, n_0+m}$. At the beginning of this copying process, $c(I_{1,1})$ and $c(I_{1,7})$ are equal to the index $n+\sum_{i=1}^{i_0} m_i$ of the last copied value ($x_n$ or $y_{m_0}$). Therefore, we need $c(I_{1,6})=n_0+m$ determined by the subprogram labeled by $\overline 1$ and thus we get $c(I_{1,1})=c(I_{1,7})=n+\sum_{i=1}^{i_0} m_i+m$. By the pseudo instructions labeled by $4^*$, the last value $y_{m_0+m}$ can be copied into a further register whose index is then stored in $I_{3,k_{{\cal N}_Q}+1}$ and also in $I_{3,k_{{\cal N}_Q}+2}$. The repeated use of ${\sf init}(Z_{3,I_{3,k_{{\cal N}_Q}+1}})$ can be carried out instead of applying the input procedure of ${\cal N}_Q\in{\sf M}_{\cal A}$. 

\begin{overview}
[Program segment ${\cal P}_{\rm init}^*$ for copying the initial values]\label{CopyInp2} 

\hfill

\nopagebreak 

\noindent \fbox{\parbox{11.8cm}{

\vspace{0.1cm}

\quad {\sf 
$2^*: \, (Z_{3,I_{3,1}+1},\ldots, Z_{3,I_{1,6}}):=(Z_{1,I_{1,1}+1},\ldots, Z_{1,I_{1,7}});$

\vspace{0.1cm}

\quad $3^*: \, I_{3,1}:=I_{1,6};$} \hfill {\scriptsize (This implies $c(I_{3,1})=n_0+m$.)}

\vspace{0.1cm}

\quad $4^*:\, I_{3,k_{{\cal N}_Q}+1}:=I_{3,1}; \, I_{3,k_{{\cal N}_Q}+2}:=I_{3,k_{{\cal N}_Q}+1};$

\hspace{0.95cm} ${\sf init}(Z_{3,I_{3,k_{{\cal N}_Q}+1}});$

\vspace{0.1cm}

\quad $5^*: \, I_{1,5}:=I_{1,5}+1;$  {\sf goto $\textcolor{blue}{1''};$}}}
\end{overview}
\noindent 
The command {\sf goto $1''$} enables the start of the simulation of $s_0$ steps of ${\cal N}_Q$ on the input given by (Inp) by using the segment ${\cal P}_{{\cal N}_Q}''$ derived from ${\cal P}_{{\cal N}_Q}$ as follows.

\begin{overview}
[Program segment ${\cal P}_{{\cal N}_Q}''$ for simulating $s_0$ steps of ${{\cal N}_Q}$]\label{SubprSimPNQ}

\hfill

\nopagebreak 

\noindent \fbox{\parbox{11.8cm}{

Let the program ${\cal P}_{{\cal N}_Q}$ of the 1-tape machine ${{\cal N}_Q}\in {\sf M}_{\cal A}$ be of the form \textcolor{blue}{$(*)$} (see p.\,\,\pageref{star}) and let ${\cal P}_{{\cal N}_Q}''$ belong to ${\cal P}_{{\cal N}_{\cal M}}$ and given by

\vspace{0.2cm}
$\quad \textcolor{blue}{1 \,''} : {\sf inst} ''_1 ;\quad \ldots \quad;\,\, (\ell _{{\cal P}_{{\cal N}_Q}}-1)\,'' : {\sf inst} ''_{ \ell _{{\cal P}_{{\cal N}_Q}}-1} ;\,\, \ell _{{\cal P}_{{\cal N}_Q}}'' : {\sf goto }\,\, \textcolor{blue}{\widetilde 2} ;$
\vspace{0.2cm}

\noindent such that each {\sf inst}$''_k$ results from {\sf instruction}$_k$ by 

\quad inserting (in front) \parskip -0.5mm 
\begin{itemize} 
\item\parskip -0.5mm 
 {\sf if $I_{1,9}=I_{1,5}$ then goto \textcolor{blue}{$\overline{2}$ }else $I_{1,9}:=I_{1,9}+1;$}\end{itemize} 
\quad and replacing

\begin{itemize} 
\item\parskip -0.5mm 
 each symbol $Z_i$ by $Z_{3,i}$, 
\item 
each symbol $I_i$ by $I_{3,i}$, 
\item
each command {\sf goto $\ell$} by {\sf goto $\ell''$}
\end{itemize} \parskip -2mm 

\quad and extending {\sf inst}$''_k$ once more as follows if {\sf instruction}$_k$ is of type (7).
}}
\nopagebreak 

\noindent \fbox{\parbox{11.8cm}{
If {\sf instruction}$_k$ is an instruction of the form \,\fbox{$I_j:=I_j+1$}, then let us extend \,{\sf inst}$''_k$\, once more by adding ${\sf init}(Z_{3,I_{3,k_{{\cal N}_Q}+1}})$ such that \, ${\sf inst}''_k$ \, is the following subprogram.

\vspace{0.1cm}

\hspace{1cm}$ ${\sf if $I_{1,9}=I_{1,5}$ then goto \textcolor{blue}{$\overline{2}$} else $I_{1,9}:=I_{1,9}+1;$}

\hspace{1cm}$I_{3,j}:=I_{3,j}+1;$ ${\sf init}(Z_{3,I_{3,k_{{\cal N}_Q}+1}})$
}}\end{overview}

By executing the subprogram ${\sf init}(Z_{3,I_{3,k_{{\cal N}_Q}+1}})$, it is ensured that any $Z$-register $Z_{3,m'}$ (with $m'=c(I_{3,k_{{\cal N}_Q}+1})+1$) on the right-hand side of the input for ${\cal N}_Q$ on the third tape of ${\cal N}_{\cal M}$ --- more precisely each of the registers $Z_{3,n_0+m+1}$, $Z_{3,n_0+m+2}, \ldots$ --- contains the last component $y_{m_0+m}$ of the input for ${\cal N}_Q$, stored in $Z_{3,n_0+m}$ at the beginning, before such a register $Z_{3,m'}$ is used.

\section*{Summary}
\addcontentsline{toc}{section}{\bf Summary}
\markboth{SUMMARY}{SUMMARY}

\begin{overview}
[Several non-determinisms]

\hfill

\nopagebreak 
\noindent \fbox{\parbox{11.8cm}{\hfill {\small Let ${\cal A}$ be any first-order structure.} 

\vspace{0.1cm}

Non-deterministic semi-decidability 

\vspace{0.2cm}

{\small

\qquad If $Q \,\,\textcolor{blue}{\subseteq}\,\, U_{\cal A}^\infty$ and $Q\in {\rm SDEC}_{\cal A}$, then 

\vspace{0.2cm}

\hspace{1.6cm}$ ({\rm SDEC}_{\cal A}^\nu)^Q\hspace{0.15cm} \,\,\textcolor{blue}{\subseteq}\,\,\hspace{0.15cm}{\rm SDEC}_{\cal A}^{\rm ND}$.
\hfill{\scriptsize (Theo.\ref{OracleMNondetM}, Rem.\ref{RemOracleMNondetM})}

\vspace{0.2cm}

\qquad If $Q \,\,\textcolor{blue}{=}\,\, U_{\cal A}^\infty$, then 

\vspace{0.2cm}

\hspace{1.6cm}$ ({\rm SDEC}_{\cal A}^\nu)^Q\hspace{0.15cm} \,\,\textcolor{blue}{=}\,\,\hspace{0.15cm}{\rm SDEC}_{\cal A}^{\rm ND}$.
\hfill{\scriptsize (Theo.\ref{OracleMEqNondM})}

\vspace{0.2cm}
}

Non-deterministic computability 

\vspace{0.2cm}

{\small

\qquad If $Q \,\,\textcolor{blue}{\subseteq}\,\, U_{\cal A}^\infty$ and $Q\in {\rm SDEC}_{\cal A}$, then 

\vspace{0.2cm}

\hspace{1.6cm}$ \{{\rm Res}_{\cal M}\mid {\cal M}\in {\sf M}_{\cal A}^\nu(Q)\}\hspace{0.15cm} \,\,\textcolor{blue}{\subseteq}\,\,\hspace{0.15cm}\{{\rm Res}_{\cal M}\mid {\cal M}\in {\sf M}_{\cal A}^{\rm ND}\}$. 
\hfill{\scriptsize (Theo.\ref{NonDetClasses1})}

\vspace{0.2cm}

\qquad If $Q \,\,\textcolor{blue}{=}\,\, U_{\cal A}^\infty$, then 

\vspace{0.2cm}

\hspace{1.6cm}$ \{{\rm Res}_{\cal M}\mid {\cal M}\in {\sf M}_{\cal A}^\nu(Q)\}\hspace{0.15cm} \,\,\textcolor{blue}{=}\,\,\hspace{0.15cm}\{{\rm Res}_{\cal M}\mid {\cal M}\in {\sf M}_{\cal A}^{\rm ND}\}$. 
\hfill{\scriptsize (Theo.\ref{NonDetClasses1})}

\vspace{0.3cm}

}}}
\end{overview}

\section*{Outlook: Non-determinisms and two constants}\label{Outlook}
\addcontentsline{toc}{section}{\bf Outlook and Acknowledgment}
\markboth{OUTLOOK}{OUTLOOK}

The power of non-deterministic machines can depend on the number of constants belonging to a structure and other basic properties. In Part II\,b, we consider in particular structures containing two different constants such as $c_1$ and $c_2$. 

\begin{overview}
[Basic properties: Identity and constants] 

\hfill

\nopagebreak 

\noindent \fbox{\parbox{11.8cm}{

\hfill {\small Let ${\cal A}$ contain the constants $c_1$ and $c_2$.}

\qquad $(a)$ \, {\small ${\cal A}$ contains the identity.}

\qquad $(b)$ \, {\small The identity is decidable over ${\cal A}$.} 

\qquad $(c)$ \, {\small The identity is semi-decidable over ${\cal A}$.} 

\qquad $(d)$ \, {\small $\{c_1\}$ and $\{c_2\}$ are semi-decidable by two machines over ${\cal A}$.}

\qquad $(e)$ \, {\small $\{c_1,c_2\}^\infty$ and thus $\{c_1,c_2\}$ are semi-decidable over ${\cal A}$.}

\vspace{0.2cm}}}

\end{overview}

\begin{overview}
[Decidability of $P$ and computability of $\chi_P$]\label{DecAndChar}

\hfill

\nopagebreak 

\noindent \fbox{\parbox{11.8cm}{

\hfill {\small Let ${\cal A}$ contain the constants $c_1$ and $c_2$ and $P\subseteq U_{\cal A}^\infty$.}

\vspace{0.2cm}

\quad $\chi_P$ is computable over ${\cal A}$ and $(d)$ $ \,\,\textcolor{blue}{\Rightarrow}\,\,$ $P\in {\rm DEC}_{\cal A}$.

\vspace{0.1cm}

\quad $P\in {\rm DEC}_{\cal A}$ $ \,\,\textcolor{blue}{\Rightarrow}\,\,$ $\chi_P$ is computable over ${\cal A}$.

\vspace{0.2cm}
}}
\end{overview}

\vspace{0.4cm}

\begin{overview}
[Several non-determinisms and the result functions]

\hfill

\noindent \fbox{\parbox{11.8cm}{

\hfill {\small Let ${\cal A}$ contain two constants $c_1$ and $c_2$.} 

\vspace{0.1cm}

Non-deterministic computability 

{\small
\qquad If $(e)$, then 

\vspace{0.1cm}

\hspace{1.6cm}$\{{\rm Res}_{\cal M}\mid {\cal M}\in {\sf M}_{\cal A}^{\rm DND}\}\hspace{0.3cm} \,\,\textcolor{blue}{\subseteq}\,\,\hspace{0.15cm} \{{\rm Res}_{\cal M}\mid {\cal M}\in {\sf M}_{\cal A}^{\rm ND}\}$. 

\vspace{0.2cm}
}}}

\noindent \fbox{\parbox{11.8cm}{

\hfill {\small Let ${\cal A}$ contain two constants $c_1$ and $c_2$ and $Q_0 \,\,\textcolor{blue}{=}\,\,\{c_1,c_2\}^2$.}

\vspace{0.1cm}

Digitally non-deterministic computability

{\small
\qquad If $(e)$, then 

\vspace{0.1cm}

\hspace{1.6cm}$\{{\rm Res}_{\cal M}\mid {\cal M}\in {\sf M}_{\cal A}^\nu(Q_0)\} \,\,\textcolor{blue}{=}\,\,\hspace{0.15cm}\{{\rm Res}_{\cal M}\mid {\cal M}\in {\sf M}_{\cal A}^{\rm DND}\}$. 

\vspace{0.2cm}}

Binary non-deterministic computability

{\small
\vspace{0.1cm}

\hspace{1.6cm}$ \{{\rm Res}_{\cal M}\mid {\cal M}\in {\sf M}_{\cal A}^\nu(Q_0)\} \,\,\textcolor{blue}{\subseteq}\,\,\hspace{0.15cm} \{{\rm Res}_{\cal M}\mid {\cal M}\in {\sf M}_{\cal A}^{\rm NDB}\}$. 

\vspace{0.2cm}}

 Digitally non-deterministic $\nu$-computability

{\small
\qquad If $(d)$ holds, then

\vspace{0.1cm}

\hspace{1.6cm}$\{{\rm Res}_{\cal M}\mid {\cal M}\in {\sf M}_{\cal A}^{\rm NDB}\} \hspace{0.3cm} \,\,\textcolor{blue}{\subseteq}\,\,\hspace{0.15cm} \{{\rm Res}_{\cal M}\mid {\cal M}\in {\sf M}_{\cal A}^\nu(Q_0)\}$. 
\vspace{0.4cm}}
}}
\end{overview}

\begin{overview}
[Semi-decidability for two types of machines]

\hfill

\nopagebreak 

\noindent \fbox{\parbox{11.8cm}{

{\small \hfill Let ${\cal A}$ be any first-order structure.

\hspace{2cm}${\rm SDEC}_{\cal A} \,\,\textcolor{blue}{=}\,\,{\rm SDEC}_{\cal A}^{\rm NDB}$.} 

\vspace{0.2cm}
}}
\end{overview}

In Part III, we will also use properties resulting from Consequences 2.2 in Part I which say that it is also possible to simulate machines over ${\cal A}_{\mathbb{N}}$ by executing index instructions of types (5) to (7). Later, we also discuss the complexity of machines, functions, and decision problems for first-order structures with identity and without identity and deterministic Moschovakis operators.

\section*{Acknowledgment}

\markboth{ACKNOWLEDGMENT}{ACKNOWLEDGMENT}

I would like to thank the participants of my lectures Theory of Abstract Computation and, in particular, Patrick Steinbrunner and Sebastian Bierba\ss{} for useful questions and the discussions. My thanks go also to all participants of meetings in Greifswald and Kloster and Pedro F. Valencia Vizcaíno. In particular, I would like to thank my co-authors Arno Pauly and Florian Steinberg for the discussions on operators related to several models of computation and Vasco Brattka, Philipp Schlicht, and Rupert H\"olzl for many interesting discussions.

{\small

\noindent For this article, we also used translators such as those of Google and DeepL and, for questions about English grammar, we used Microsoft Copilot.}

\vspace{4cm}

\vfill

\hfill {\scriptsize\sf C $\cdot $ H $\cdot$ T}
\end{document}